\documentclass[a4paper]{amsart}                

\usepackage[latin1]{inputenc}                  
\usepackage[T1]{fontenc}                       
\usepackage[spanish,english]{babel}            
\usepackage{graphicx}                  
\usepackage{amsmath,amssymb,amsthm}            
\usepackage{latexsym}                          
\usepackage{delarray}                          
\usepackage{bbm}                               
\usepackage{hyperref}                          
\usepackage{calrsfs,eufrak,mathrsfs,upgreek}   
\usepackage{bbding,trfsigns}                   
\usepackage{wasysym}                           
\usepackage{datetime}                          
\usepackage{color}                             

\DeclareMathAlphabet{\mathpzc}{OT1}{pzc}{m}{it} 

\newtheorem{teo}{Theorem}

\newtheorem{propo}{Prop.}
\newtheorem{remark}{Remark}

\newtheorem*{A}{Theorem A}

\author[J. Betancor]{J.J. Betancor}
\address{Departamento de An\'{a}lisis Matem\'{a}tico\\
Universidad de la Laguna\\
Campus de Anchieta, Avda. Astrof\'{\i}sico Francisco S\'{a}nchez, s/n\\
38271 La Laguna (Sta. Cruz de Tenerife), Espa\~na}
\email{jbetanco@ull.es, jcfarina@ull.es, asgarcia@ull.es}

\author[J.C. Fari\~{n}a]{J.C. Fari\~{n}a}

\author[A. Sanabria]{A. Sanabria}

\thanks{This paper is partially supported by MTM2010/17974.}
\date{\today}

\begin{document}

\title[Vector valued multivariate spectral multipliers]
{Vector valued multivariate spectral multipliers, Littlewood-Paley functions, and Sobolev spaces in the Hermite setting}



\begin{abstract}
In this paper we find new equivalent norms in $L^p(\mathbb{R}^n,\mathbb{B})$ by using multivariate Littlewood-Paley functions associated with Poisson semigroup for the Hermite operator, provided that $\mathbb{B}$ is a UMD Banach space with the property ($\alpha$). We make use of $\gamma$-radonifying operators to get new equivalent norms that allow us to obtain $L^p(\mathbb{R}^n,\mathbb{B})$-boundedness properties for (vector valued) multivariate spectral multipliers for Hermite operators. As application of this Hermite multiplier theorem we prove that the Banach valued Hermite Sobolev and potential spaces coincide.
\end{abstract}

\maketitle

\section{Introduction}

The Hermite operator (also called harmonic oscillator) $\tilde H$ on $\mathbb{R}^n$ is defined by
$$
\tilde H = -\Delta +|x|^2,
$$
where $\Delta$ denotes the usual Laplacian operator. For every $m \in \mathbb{N}$ we denote by $h_m$ the $m$-th Hermite function given by
$$
h_m(u) = (2^m m! \sqrt{\pi})^{-\frac12} P_m(u) e^{-\frac{u^2}{2}},\;\; u \in \mathbb{R},
$$
where  $P_m$ represents the $m$-th Hermite polynomial (\cite[p. 104]{Sz}).

If $k=(k_1,\dots,k_n) \in \mathbb{N}^n$, the $k$-th Hermite function $h_k$ is defined by
$$
h_k (x) = \prod_{j=1}^n h_{k_j} (x_j),\;\; x=(x_1,\dots,x_n)\in\mathbb{R}^n.
$$
We have that, for every $k=(k_1,\dots,k_n) \in \mathbb{N}^n$,
$$
\tilde H h_k = (2|k|+n) h_k,
$$
where $|k|=k_1+\cdots+k_n$. The sequence $\{h_k\}_{k\in\mathbb{N}^n}$ is an orthonormal basis in $L^2(\mathbb{R}^n)$. Moreover, the linear space $\textup{span}\{h_k\}_{k\in\mathbb{N}^n}$ generated by $\{h_k\}_{k\in\mathbb{N}^n}$ is dense in $L^p(\mathbb{R}^n)$, $1\le p<\infty$ (\cite[Lemma 2.3]{StemTo1}).

The Hermite operator $H$ is defined by
$$
Hf = \sum_{k\in\mathbb{N}^n} (2|k|+n) c_k(f) h_k,\;\; f\in D(H),
$$
where
$$
D(H) = \{f\in L^2(\mathbb{R}^n) \colon \sum_{k\in\mathbb{N}^n}(2|k|+n)^2 |c_k(f)|^2< \infty\}
$$
and, for every $k\in\mathbb{N}^n$
$$
c_k(f) = \int_{\mathbb{R}^n} h_k(x) f(x) dx,\;\; f\in L^2(\mathbb{R}^n).
$$
According to \cite[Lemma 1.2]{StemTo1}, the space $C_c^\infty (\mathbb{R}^n)$ of smooth compactly supported functions in $\mathbb{R}^n$ is contained in the domain $D(H)$ of $H$ and $Hf=\tilde H f$, $f\in C_c^\infty (\mathbb{R}^n)$.

Harmonic analysis associated with Hermite polynomial expansions were begun by Muckenhoupt (\cite{Mu2} and \cite{Mu1}) who considered the one dimensional setting. Later, Sj\"ogren (\cite{Sj}), Fabes, Guti\'errez and Scotto (\cite{FGS}) and Urbina (\cite{U}) studied harmonic analysis operators associated with the Ornstein-Uhlenbech operator in $\mathbb{R}^n$. In the last decade this topic has been studied by a host of authors (\cite{CMM}, \cite{FHS}, \cite{GCMMST}, \cite{LY} and \cite{MS}).

In the monography of Thangavelu (\cite{Th}) and in the papers of K. Stempak and J.L. Torrea (\cite{StemTo1}, \cite{StemTo2} and \cite{StemTo3}) the harmonic analysis operators associated with the Hermite operator (Hermite function expansions setting) were investigated. Later, this study have been completed in a series of papers (see \cite{AT}, \cite{BCCFR}, \cite{BFRST}, \cite{H} and \cite{L} among others).

The heat semigroup $\{W_t^H\}_{t>0}$ generated by $-H$ in $L^2(\mathbb{R}^n)$ is defined by
$$
W_t^H (f) = \sum_{k\in\mathbb{N}^n} e^{-(2|k|+n)t} c_k(f) h_k,\;\; f \in L^2(\mathbb{R}^n)\;\textup{and}\;t>0.
$$
According to Mehler's formula (\cite[Lemma 1.1.1, p. 2]{Th}) we can write, for every $t>0$ and $f\in L^2(\mathbb{R}^n)$,
\begin{equation}\label{1.1}
W_t^H (f) (x) = \int_{\mathbb{R}^n} W_t^H(x,y) f(y) dy,
\end{equation}
where
$$
W_t^H(x,y) = \frac{1}{\pi^\frac{n}{2}} \left( \frac{e^{-2t}}{1-e^{-4t}}\right)^{\frac{n}{2}} e^{ -\frac14\left( |x-y|^2 \frac{1+e^{-2t}}{1-e^{-2t}} + |x+y|^2 \frac{1-e^{-2t}}{1+e^{-2t}}\right)}, \;\; x,y\in\mathbb{R}^n\;\textup{and}\;t>0.
$$
Moreover, $\{W_t^H\}_{t>0}$ (defined by (\ref{1.1})) is the semigroup of contractions in $L^p(\mathbb{R}^n)$, $1\le p<\infty$, generated by $-H$.

The Poisson semigroup $\{P_t^H\}_{t>0}$ associated with the Hermite operator (generated by $-\sqrt{H})$ in $L^p(\mathbb{R}^n)$, $1\le p<\infty$, is defined by using the subordination formula as follows; for every $t>0$,
$$
P_t^H(f)(x) = \frac{t}{2\sqrt{\pi}} \int_0^\infty \frac{e^{-\frac{t^2}{4u}}}{u^{\frac{3}{2}}} W_u^H(f)(x) du,\;\; f \in L^p(\mathbb{R}^n).
$$
$\{P_t^H\}_{t>0}$ is a semigroup of contractions in $L^p(\mathbb{R}^n)$, $1\le p<\infty$.

The Hermite semigroups $\{W_t^H\}_{t>0}$ and $\{P_t^H\}_{t>0}$ are not conservative. Then, $\{W_t^H\}_{t>0}$ and $\{P_t^H\}_{t>0}$ are not diffusion semigroups (in the sense of Stein \cite{St}).

The square function (also called Littlewood-Paley function) $g_{W,k}^H$ associated with the Hermite semigroup $\{W_t^H\}_{t>0}$ is defined by
$$
g_{W,k}^H (f)(x) = \left( \int_0^\infty |t^k \partial_t^k W_t^H(f)(x) |^2 \frac{dt}{t}\right)^{\frac12},\;\; x \in\mathbb{R}^n,
$$
for every $f\in L^p(\mathbb{R}^n)$ and $k\in\mathbb{N}\setminus\{0\}$. $L^p$-boundedness properties of the (nonlinear, but almost linear) operator $g_{W,k}^H$, $k\in\mathbb{N}\setminus\{0\}$, were established in \cite[Theorem 2.3.2, p. 41]{Th} and \cite[Theorem 2.2]{StemTo2}. By using $g_{W,k}^H$ it is possible to define equivalent norms in $L^p(\mathbb{R}^n)$, $1<p<\infty$. Indeed, for every $1<p<\infty$ and $k\in\mathbb{N}\setminus\{0\}$, there exists $C>0$ such that
\begin{equation}\label{1.2}
    \frac{1}{C} \Vert f \Vert_{L^p(\mathbb{R}^n)} \le \Vert g_{W,k}^H (f) \Vert_{L^p(\mathbb{R}^n)}\le C \Vert f \Vert_{L^p(\mathbb{R}^n)},\;\; f\in L^p(\mathbb{R}^n).
\end{equation}
Equivalence (\ref{1.2}) allows to obtain $L^p$-boundedness for spectral Hermite multipliers (see \cite[Theorem 2.4.1, p. 45]{Th}).

The Littlewood-Paley functions associated with the Hermite-Poisson semigroup $\{P_t^H\}_{t>0}$ are defined by
$$
g_{P,k}^H (f)(x) = \left( \int_0^\infty |t^k \partial_t^k P_t^H(f)(x) |^2 \frac{dt}{t}\right)^{\frac12},\;\; x \in\mathbb{R}^n,
$$
for every $f\in L^p(\mathbb{R}^n)$, $1<p<\infty$, and $k\in\mathbb{N}\setminus\{0\}$. Equivalence (\ref{1.2}) also holds when $g_{W,k}^H$ is replaced by $g_{P,k}^H$, for every $1<p<\infty$ and $k\in\mathbb{N}\setminus\{0\}$.

In this paper we are interested in multivariate spectral Hermite multipliers. We need to consider multiparameter square functions associated with the Hermite operator. Littlewood-Paley functions of this kind have been studied recently by Wr\'obel (\cite{Wr} and \cite{Wr1}). Although our results can also be established by using square functions for the Hermite heat semigroup, the manipulations and calculations are simpler when we consider square functions for the Hermite-Poisson semigroup. Then, we will use Littlewood-Paley functions associated with the Hermite-Poisson semigroup in the sequel.

Let $k=(k_1,\dots,k_n)\in (\mathbb{N}\setminus\{0\})^n$. We define the multivariate square function $g_{P,k}^H$ as follows
$$
g_{P,k}^H (f)(x) = \left( \int_{(0,\infty)^n} \left|\int_{\mathbb{R}^n} \prod_{j=1}^n t_j^{k_j} \partial_{t_j}^{k_j} P_{t_j}^H (x_j,y_j) f(y_1,\dots,y_n) dy \right|^2 \frac{dt_1\cdots dt_n}{t_1\cdots t_n} \right)^{\frac12},\;\; x\in \mathbb{R}^n,
$$
for every $f\in L^p(\mathbb{R}^n)$, $1<p<\infty$.

By proceeding as in the proof of \cite[Theorem 2.4]{Wr} we can show the following property.
\begin{A}
Let $1<p<\infty$ and $k=(k_1,\dots,k_n) \in (\mathbb{N}\setminus\{0\})^n$. Then, there exists $C>0$ such that
\begin{equation}\label{1.3}
\frac{1}{C}||f||_{L^p(\mathbb{R}^n)} \le ||g^H_{P,k}(f)||_{L^p(\mathbb{R}^n)} \le C ||f||_{L^p(\mathbb{R}^n)},\;\;
f\in L^p(\mathbb{R}^n).
\end{equation}
\end{A}

\bigskip
Equivalence in (\ref{1.3}) is the key to establish a multivariate spectral multiplier theorem (\cite[Theorem 2.2]{Wr}).

Our first objective in this paper is to prove a Banach valued version of Theorem A. Then, inspired by the results of Meda \cite{Me} and Wr\'obel (\cite{Wr} and \cite{Wr1}), we prove a vector valued multivariate spectral multiplier theorem in the Hermite setting. As an application of our Hermite multiplier theorem we see that the Banach valued Hermite Sobolev and Hermite potential spaces coincide provided that the Banach space has UMD and ($\alpha$) properties.

Let $\mathbb{B}$ be a Banach space and $k=(k_1,\dots,k_n)\in (\mathbb{N}\setminus\{0\})^n$. In order to define Hermite square functions in the $\mathbb{B}$-valued setting, the more natural way is to replace the absolute value that appears in the scalar case by the norm in $\mathbb{B}$. With this idea we define, for every $f\in L^p(\mathbb{R}^n,\mathbb{B})$, $1<p<\infty$,
$$
g_{P,k;\mathbb{B}}^H (f)(x) = \left( \int_{(0,\infty)^n} \left\Vert \int_{\mathbb{R}^n} \prod_{j=1}^n t_j^{k_j} \partial_{t_j}^{k_j} P_{t_j}^H (x_j,y_j) f(y_1,\dots,y_n) dy \right\Vert_{\mathbb{B}}^2 \frac{dt_1\cdots dt_n}{t_1\cdots t_n} \right)^{\frac12},\;\; x\in \mathbb{R}^n.
$$
This kind of Littlewood-Paley function have been considered in \cite{BFRST2}, \cite{HTV}, \cite{MTX}, \cite{TZ} and \cite{Xu}, among others, always in the univariate case. According to the results in \cite{Kwa} and \cite{MTX}, even in the case $n=k=1$, there exists $C>0$ such that, for some $1<p<\infty$,
$$
\frac{1}{C} \Vert f \Vert_{L^p(\mathbb{R},\mathbb{B})} \le \Vert g_{P,1;\mathbb{B}}^H (f) \Vert_{L^p(\mathbb{R}^n)} \le C \Vert f \Vert_{L^p(\mathbb{R},\mathbb{B})},\;\; f\in L^p(\mathbb{R},\mathbb{B}),
$$
if, and only if, $\mathbb{B}$ is isomorphic to a Hilbert space.

If we want to get equivalent norms in $L^p(\mathbb{R}^n,\mathbb{B})$ by using multivariate Hermite square functions for Banach spaces $\mathbb{B}$ that are not isomorphic to Hilbert spaces, we need to follow other way. This question in the univariate case have been studied by Hyt\"onen \cite{HyRMI} and Kaiser and Weis \cite{KW}. Hyt\"onen in \cite{HyRMI} used stochastic integration to define Littlewood-Paley functions associated with diffusion semigroups. We recall that the Hermite heat and Poisson semigroups are not diffusion semigroups. Kaiser and Weis \cite{KW} studied vector valued square functions defined by convolutions by using $\gamma$-radonifying operators. The two procedure are, in some sense, equivalent (see \cite{VanNee}). The results obtained in \cite{HyRMI} and \cite{KW} are valid for UMD Banach spaces. It is well-known that all the Banach spaces that are isomorphic to Hilbert spaces are also UMD but there exists UMD Banach spaces that are not isomorphic to Hilbert spaces (for instance, $L^p(\mathbb{R})$, $1<p<\infty$, $p\ne 2$).

In this paper we define multivariate square functions associated to Hermite operators by using $\gamma$-radonifying operators following the ideas in \cite{KW}.

We now recall some definitions and properties that will be useful in the sequel. We consider the Hilbert space $\mathcal{H}^n = L^2\left((0,\infty)^n,\frac{dt_1\cdots dt_n}{t_1\cdots t_n}\right)$, $n\in\mathbb{N}\setminus\{0\}$,  and we choose an orthonormal basis $\{\varphi_j\}_{j\in\mathbb{N}}$ of $\mathcal{H}^n$. Suppose that $\{\gamma_j\}_{j\in \mathbb{N}}$ is a sequence of independent standard Gaussians on a probability space $(\Omega,\mathbb{P})$ and that $\mathbb{B}$ is a Banach space. By $\mathcal{L}(\mathcal{H}^n,\mathbb{B})$ we denote the space of bounded linear operator from $\mathcal{H}^n$ into $\mathbb{B}$. We say that $T\in \mathcal{L}(\mathcal{H}^n,\mathbb{B})$ is a $\gamma$-radonifying operator (written $T\in\gamma(\mathcal{H}^n,\mathbb{B})$) when the series $\sum_{j=1}^\infty \gamma_j T(\varphi_j)$ converges in $L^2(\Omega,\mathbb{B})$.

This definition does not depend on the sequences $\{\varphi_j\}_{j\in \mathbb{N}}$ and $\{\gamma_j\}_{j\in \mathbb{N}}$. By endowing $\gamma(\mathcal{H}^n,\mathbb{B})$ with the norm $\Vert \cdot \Vert_{\gamma(\mathcal{H}^n,\mathbb{B})}$ defined by
$$
\Vert T \Vert_{\gamma(\mathcal{H}^n,\mathbb{B})} =  \left( \mathbb{E} \Big\Vert \sum_{j=1}^\infty \gamma_j T(\varphi_j) \Big \Vert_{\mathbb{B}}^2 \right)^{\frac12},\;\;T\in \gamma(\mathcal{H}^n,\mathbb{B}),
$$
which is independent of the basis $\{\varphi_j\}_{j=1}^\infty$, $\gamma(\mathcal{H}^n,\mathbb{B})$ becomes a Banach space. Also, we have that
$$
\Vert T \Vert_{\gamma(\mathcal{H}^n,\mathbb{B})} = \sup \left( \mathbb{E} \Big \Vert \sum_{j=1}^l \gamma_j T(\varphi_j) \Big \Vert_{\mathbb{B}}^2 \right)^{\frac12},\;\;T\in \gamma(\mathcal{H}^n,\mathbb{B}),
$$
where the supremum is taken over all the finite orthonormal sets $\{\varphi_j\}_{j=1}^l$ in $\mathcal{H}^n$.  Moreover, if the Banach space $\mathbb B$ does not contain any copies of $c_0$ (for instance, when $\mathbb B$ is a UMD space), then $T\in \mathcal{L}(\mathcal{H}^n,\mathbb B)$ if and only if
$$
\sup\left(\mathbb E\Big\Vert\sum_{j=1}^l \gamma_j T(\varphi_j)\Big\Vert_{\mathbb B}^2\right)^{\frac{1}{2}} < \infty,
$$
where the supremum is taken over all the finite orthonormal sets $\{\varphi_j\}_{j=1}^l$ in $\mathcal{H}^n$.

Suppose that $f\colon (0,\infty)^n \to \mathbb{B}$ is a strongly measurable function. If $f$ is weakly $\mathcal{H}^n$, that is, for every $S\in \mathbb{B}^*$, $S\circ f \in \mathcal{H}^n$, then there exists $T_f\in \mathcal{L}(\mathcal{H}^n,\mathbb{B})$ satisfying that
$$
\langle S, T_f(\varphi) \rangle_{\mathbb{B}^*,\mathbb{B}} = \int_{(0,\infty)^n} \langle S, f(t) \rangle_{\mathbb{B}^*,\mathbb{B}} \varphi(t) \frac{dt_1\cdots dt_n}{t_1\cdots t_n},\;\; \varphi \in \mathcal{H}^n.
$$
We say that $f\in\gamma\Big((0,\infty)^n,\frac{dt_1\cdots dt_n}{t_1\cdots t_n},\mathbb{B}\Big)$ when $T_f \in \gamma(\mathcal{H}^n,\mathbb{B})$. By identifying $f$ with $T_f$, for every $f\in \gamma\Big((0,\infty)^n,\frac{dt_1\cdots dt_n}{t_1\cdots t_n},\mathbb{B}\Big)$, $\gamma\Big((0,\infty)^n,\frac{dt_1\cdots dt_n}{t_1\cdots t_n},\mathbb{B}\Big)$ is a dense subspace of $\gamma(\mathcal{H}^n,\mathbb{B})$ provided that $\mathbb{B}$ does not contain any copy of $c_0$. If $T_f \not\in \gamma(\mathcal{H}^n,\mathbb{B})$, or if $f$ is not even weakly $\mathcal{H}^n$, then we write $\Vert f \Vert_{\gamma(\mathcal{H}^n,\mathbb{B})} = \infty$. The main properties of $\gamma$-radonifying operators can be found in \cite{VanNee}.

As it is well-known the Hilbert transform is defined by
$$
\mathrm{H}(f)(x) = \lim_{\epsilon \to 0} \frac{1}{\pi} \int_{|x-y|>\epsilon} \frac{f(y)}{y-x} dy,\;\;\textup{a.e.}\; x \in \mathbb{R},
$$
for every $f\in L^p(\mathbb{R})$, $1\le p<\infty$, and it is a bounded operator from $L^p(\mathbb{R})$ into itself, for every $1<p<\infty$, and from $L^1(\mathbb{R})$ into $L^{1,\infty}(\mathbb{R})$. The operator $\mathrm{H}\otimes \textup{Id}_{\mathbb{B}}$ is defined as usual in $L^p(\mathbb{R}) \otimes \mathbb{B}$. The Banach space $\mathbb{B}$ is said to be UMD when the operator $\mathrm{H}\otimes \textup{Id}_{\mathbb{B}}$ can be extended to $L^p(\mathbb{R},\mathbb{B})$ as a bounded operator from $L^p(\mathbb{R},\mathbb{B})$ into itself for some (equivalently, for every $1<p<\infty$). Main results about UMD Banach spaces were established by Bourgain (\cite{Bou}), Burkholder (\cite{Bu}) and Rubio de Francia (\cite{Ru}).

If $\{\epsilon_j\}_{j=1}^\infty$ is a sequence of independent symmetric $\pm 1$-valued random variables (usually called Rademacher variables) on some probability space, we denote by $\mathbb{E}_\epsilon$ the corresponding expectation operator.

Suppose that $\{\epsilon_j\}_{j=1}^\infty$ and $\{\eta_j\}_{j=1}^\infty$ are two independent sequences of Rademacher variables. We say that a Banach space $\mathbb{B}$ has (Pisier's) property ($\alpha$) when there exists $C>0$ such that
$$
\mathbb{E}_\epsilon \mathbb{E}_\eta \Big\Vert \sum_{i,j=1}^N \alpha_{i,j} \epsilon_i \eta_j x_{i,j} \Big\Vert_{\mathbb{B}} \le C\mathbb{E}_\epsilon \mathbb{E}_\eta  \Big\Vert \sum_{i,j=1}^N \epsilon_i \eta_j x_{i,j} \Big\Vert_{\mathbb{B}},
$$
for every $\alpha_{i,j}\in\{+1,-1\}$, $x_{i,j} \in \mathbb{B}$, $i,j=1,\dots, N$, and $N \in \mathbb{N}\setminus\{0\}$. This property is satisfied by the commutative $L^p$ spaces, $1\leq p <\infty$, and it is also inherited from $X$ to $L^p(X)$, $1 \leq p < \infty$. A Banach lattice has the property $(\alpha)$ if and only if it has finite cotype.

UMD and ($\alpha$) properties of Banach spaces are crucial in order to prove Banach valued Fourier multipliers theorems of Mikhlin type (see \cite{HW} and \cite{We}, among others)
and joint $H^\infty$ functional calculus (\cite{La}).

Segovia and Wheeden (\cite[p. 248]{SW}) introduced the notion of fractional derivative $\partial_t^\alpha,~\alpha>0$, as follows. Suppose that $\alpha >0$ and $m\in\mathbb{N}\setminus\{0\}$ such that $m-1\le \alpha <m$. If $f$ is a reasonable nice function on $(0,\infty)\times \mathbb{R}^n$ we define
$$
\partial_t^\alpha f(t,x) = \frac{e^{-\pi (m-\alpha) i}}{\Gamma(m-\alpha)} \int_0^\infty \partial_t^m f(t+s,x) s^{m-\alpha-1} ds,\;\; x \in \mathbb{R}^n~\textup{and}~t\in (0,\infty).
$$
Littlewood-Paley functions involving fractional derivatives were used in \cite{SW} to characterize classical Sobolev spaces. Fractional square functions associated with diffusion semigroups have been studied in \cite{AST} and \cite{TZ} in a Banach valued setting.

We define the univariate fractional $g$-function operator $\mathcal{G}_{P,\alpha;\mathbb{B}}^H$, $\alpha>0$, associated with the Hermite-Poisson semigroup by
$$
\mathcal{G}_{P,\alpha;\mathbb{B}}^H (f)(t,x) = t^\alpha \partial_t^\alpha P_t^H (f)(x),\;\; x\in \mathbb{R}^n~\textup{and}~t >0,
$$
for every $f\in L^p(\mathbb{R}^n,\mathbb{B})$.

In \cite[Theorem 1]{BCCFR1} it was established that $\mathcal{G}_{P,\alpha;\mathbb{B}}^H$ allows us to get equivalent norms in $L^p(\mathbb{R}^n,\mathbb{B})$, $1<p<\infty$.

\begin{teo}\label{nteo1}
Let $\mathbb{B}$ be a UMD Banach space. For every $\alpha >0$ and $1<p<\infty$ there exists $C>0$ such that
$$
\frac{1}{C} \Vert f \Vert_{L^p(\mathbb{R}^n,\mathbb{B})} \le \Vert \mathcal{G}_{P,\alpha;\mathbb{B}}^H (f) \Vert_{L^p(\mathbb{R}^n,\gamma(\mathcal{H}^1,\mathbb{B}))} \le C \Vert f \Vert_{L^p(\mathbb{R}^n,\mathbb{B})},\;\; f \in L^p(\mathbb{R}^n,\mathbb{B}).
$$
\end{teo}

\bigskip
Our first result is a multivariate version of Theorem \ref{nteo1}. If $k=(k_1,\dots,k_n) \in (\mathbb{N}\setminus\{0\})^n$, we consider the $g$-function associated with the Hermite operator defined by
$$
G_{P,k;\mathbb{B}}^H(f)(t,x) = \int_{\mathbb{R}^n} \prod_{j=1}^n t_j^{k_j} \partial_{t_j}^{k_j} P_{t_j}^H (x_j,y_j) f(y_1,\dots,y_n) dy_1\dots dy_n,
$$
for every $f\in L^p(\mathbb{R}^n,\mathbb{B})$, $1<p<\infty$, where $t=(t_1,\dots,t_n)\in(0,\infty)^n$ and $x=(x_1,\dots,x_n)\in\mathbb{R}^n$.

\begin{teo}\label{theorem1.4}
Let $\mathbb{B}$ be a UMD Banach space with the property ($\alpha$), $k\in (\mathbb{N}\setminus\{0\})^n$ and $1<p<\infty$. Then, there exists $C>0$ such that
$$
\frac{1}{C}||f||_{L^p(\mathbb{R}^n, \mathbb{B})} \le ||G_{P,k;\mathbb{B}}^H (f)||_{L^p(\mathbb{R}^n, \gamma(\mathcal{H}^n,\mathbb{B}))} \le C ||f||_{L^p(\mathbb{R}^n, \mathbb{B})},\;\;
f\in L^p(\mathbb{R}^n, \mathbb{B}).
$$
\end{teo}

\bigskip
Since $\gamma(\mathcal{H}^n,\mathbb{C}) = \mathcal{H}^n$, Theorem A can be seen as a special (scalar) case of Theorem \ref{theorem1.4}.

In order to prove Theorem \ref{theorem1.4} we use that $\gamma(\mathcal{H}^n,\mathbb{B}) \simeq \gamma(\mathcal{H}^l,\gamma(\mathcal{H}^{n-l},\mathbb{B}))$, for every $l\in\mathbb{N}$, $1\le l\le n-1$ (\cite[Corollary 3.5, (a)]{vNW}). In some sense, this property characterizes the property ($\alpha$) for the Banach space $\mathbb{B}$ (see \cite[Corollary 3.5, (2a) and (2b)]{vNW}).

\begin{remark}
Equivalence of norms for $L^p(\mathbb{R}^n,\mathbb{B})$ established in Theorem \ref{theorem1.4} also holds when we consider the more general multivariate $g$-function operator involving fractional derivatives defined as follows. Suppose that $n_j \in\mathbb{N}\setminus\{0\}$, $j=1,\dots,l \in \mathbb{N}$, such that $\sum_{j=1}^l n_j = n$, and $\alpha =(\alpha_1,\dots,\alpha_l)$ being $\alpha_j>0$, $j=1,\dots,l$. We define the Littlewood-Paley type operator by
$$
G_{P,\alpha;\mathbb{B}}^H (f)(t,x) = \int_{\mathbb{R}^n} \prod_{j=1}^l t_j^{\alpha_j} \partial_{t_j}^{\alpha_j} P_{t_j}^H (x^j,y^j) f(y^1,\dots,y^l) dy,
$$
where $t=(t_1,\dots,t_l)\in(0,\infty)^l$ and $x=(x^1,\dots,x^l)\in\mathbb{R}^{n_1}\times \dots \times \mathbb{R}^{n_l}$.

We prefer to state Theorem \ref{theorem1.4} in the present form because this simpler one is sufficient to obtain our results about vector valued Hermite Sobolev spaces (see Theorem \ref{theorem1.7} bellow). Moreover, the proof of the corresponding equivalence of norms for the more general operator $G_{P,\alpha;\mathbb{B}}^H$, $\alpha \in (0,\infty)^l$, can be made by using the same ideas used to show Theorem \ref{theorem1.4}, with straightforward notational elements and manipulations.
\end{remark}

Suppose now that $\textrm{m}$ is a bounded Borel measurable function from $(0,\infty)^n$ into $\mathbb{C}$. The Hermite multivariate multiplier $T_\textrm{m}$ associated with $\textrm{m}$ is defined by
$$
T_\textrm{m} f = \sum_{k=(k_1,\dots,k_n)\in\mathbb{N}^n} \textrm{m}(\lambda_{k_1},\cdots,\lambda_{k_n}) c_k(f) h_k,\;\; f\in L^2(\mathbb{R}^n),
$$
where $\lambda_l=2l+1$, for every $l\in\mathbb{N}$. $T_\textrm{m}$ is a bounded operator from $L^2(\mathbb{R}^n)$ into itself. Thangavelu \cite[Theorem 4.2.1]{Th} established Mikhlin-H\"ormander type conditions on $\textrm{m}$ under that $T_\textrm{m}$ is bounded from $L^p(\mathbb{R}^n)$ into itself, for every $1<p<\infty$. We define the operator $T_\textrm{m} \otimes \textup{Id}_{\mathbb{B}}$ on $L^2(\mathbb{R}^n) \otimes \mathbb{B}$ in the usual way. Motivated by the results of Meda \cite{Me} and Wr\'obel (\cite{Wr} and \cite{Wr1}) we establish the following Hermite multivariate multiplier result in a Banach valued setting. Previously we introduce some notations (see \cite{Me}). For every $\alpha = (\alpha_1,\dots,\alpha_n) \in \mathbb{N}^n$, we consider
$$
m_\alpha (t,\lambda) = \prod_{j=1}^n (t_j \lambda_j)^{\alpha_j} e^{-\frac{t_j \lambda_j}{2}} M(\lambda),
$$
where $t=(t_1,\dots,t_n)\in (0,\infty)^n$, $\lambda=(\lambda_1,\dots,\lambda_n) \in \mathbb{R}^n$ and $M(\lambda)=\textrm{m}(\lambda_1^2,\dots,\lambda_n^2)$, $\lambda=(\lambda_1,\dots,\lambda_n)\in\mathbb{R}^n$. We define
$$
\mathcal{M}_\alpha (t,u) = \int_{(0,\infty)^n} \prod_{j=1}^n \lambda_j^{-i u_j -1} m_\alpha (t,\lambda) d\lambda,
$$
with $t=(t_1,\dots,t_n)\in (0,\infty)^n$ and $u=(u_1,\dots,u_n)\in\mathbb{R}^n$.

Let $\beta=(\beta_1,\dots,\beta_n)\in \mathbb{N}^n$. We write $L^{i\beta}$ to refer us to the operator $T_{\textrm{m}_\beta}$ where $\textrm{m}_\beta (\lambda_1,\dots,\lambda_n) = \prod_{j=1}^n \lambda_j^{i\beta}$, $\lambda=(\lambda_1,\dots,\lambda_n) \in (0,\infty)^n$, that is,
$$
L^{i \beta} f = \sum_{k_1,\dots,k_n=0}^\infty \prod_{j=1}^n \lambda_{k_j}^{i\beta_j} c_k (f) h_k, \;\; f\in L^2(\mathbb{R}^n).
$$
The operator $L^{i \beta} \otimes \textup{Id}_{\mathbb{B}}$ is defined in the usual way on $L^2(\mathbb{R}^n) \otimes \mathbb{B}$. By using \cite[Theorem 1.2]{BCCR} we can see that $L^{i \beta} \otimes \textup{Id}_{\mathbb{B}}$ can be extended to $L^p(\mathbb{R}^n,\mathbb{B})$ as a bounded operator from $L^p(\mathbb{R}^n,\mathbb{B})$ into itself, for every $1<p<\infty$, provided that $\mathbb{B}$ is a UMD Banach space.
\begin{teo}\label{theorem1.5}
Let $\mathbb{B}$ be a UMD Banach space with the property ($\alpha$) and $1<p<\infty$. Suppose that $\mathrm{m}$ is a bounded Borel measurable function on $(0,\infty)^n$, such that for some $\gamma\in\mathbb{N}^n$,
\begin{equation}\label{1.6}
\int_{\mathbb{R}^n} \sup_{t \in (0,\infty)^n} |\mathcal{M}_\gamma (t,u)| \Vert L^{i\frac{u}{2}} \Vert_{L^p(\mathbb{R}^n,\mathbb{B}) \to L^p(\mathbb{R}^n,\mathbb{B})} du <\infty.
\end{equation}
Then, the multiplier operator $T_\mathrm{m}$ is bounded from $L^p(\mathbb{R}^n,\mathbb{B})$ into itself.
\end{teo}

We consider a class of UMD Banach spaces, called intermediate UMD spaces, that includes all known examples of UMD spaces. We say that $\mathbb{B}$ is an intermediate UMD space when $\mathbb{B}$ is isomorphic to a closed subquotient of a complex interpolation space $[X,\mathcal{Q}]_\theta$, where $\theta \in (0,1)$, $X$ is a UMD Banach space and $\mathcal{Q}$ is a Hilbert space. This class of UMD spaces has been used recently by Berkson and Gillespie \cite{BG}, Hyt\"onen \cite{HyRMI}, Hyt\"onen and Lacey \cite{HL} and Taggart \cite[Theorem 1.5]{Ta1}. All UMD lattices are intermediate UMD spaces (\cite[Corollary on p. 216]{Ru}). It is, as far as we know, an open problem if every UMD space is an intermediate UMD space. This question was posed by Rubio de Francia \cite{Ru}. We need to use intermediate UMD spaces in the following theorem to get a suitable estimate for the operator norm $\Vert L^{i\gamma}\Vert_{L^p(\mathbb{R}^n,\mathbb{B}) \to L^p(\mathbb{R}^n,\mathbb{B})}$, $\gamma \in \mathbb{N}^n$ and $1<p<\infty$.

For every $\psi \in (0,\pi)$ we denote by $\Gamma_\psi$ the $\mathbb{C}^n$-sector
$$
\Gamma_\psi = \{(z_1,\dots,z_n)\in \mathbb{C}^n \colon |\textup{Arg}\;z_j|<\psi,\; j=1,\dots,n\}.
$$
In the following theorem we specify conditions in order that a function $\textrm{m}$ satisfies the property (\ref{1.6}) in Theorem \ref{theorem1.5}.
\begin{teo}\label{theorem1.6}
Suppose that $\mathbb{B}$ is isomorphic to a closed subquotient of $[X,\mathcal{Q}]_\theta$, where $\theta\in (0,1)$, $X$ is a UMD Banach space and $\mathcal{Q}$ is a Hilbert space and that $\mathbb{B}$ has the property ($\alpha$). If $\mathrm{m}$ is a bounded holomorphic function in $\Gamma_\psi$ for some $\psi > \frac{\pi}{4}$, then the multiplier operator $T_\mathrm{m}$ can be extended to $L^p(\mathbb{R}^n,\mathbb{B})$ as a bounded operator from $L^p(\mathbb{R}^n,\mathbb{B})$ into itself, provided that $\left|\frac{2}{p}-1\right|<\theta$.
\end{teo}

The arguments used in the proofs of Theorems \ref{theorem1.5} and \ref{theorem1.6} can also be employed to show vector valued versions of \cite[Theorems 1 and 3]{Me} for the Hermite operator. In this case we must apply \cite[Theorem 4.2]{KW} and it is sufficient to consider UMD Banach spaces (not necessarily with the property ($\alpha$)).  We point out that the fundamental condition $\left|\frac{2}{p}-1\right|< \theta$ in Theorem \ref{theorem1.6} comes from \cite[Corollary 2.5.3]{TA} where the $L^p$-norm of imaginary powers of operators is estimated in a Banach valued setting. In the scalar case, by tacking into account the behavior of the $L^p$-norm of the imaginary power of our Hermite operator (polynomial growth, see \cite[(4.3)]{Wr} for instance), if $m$ satisfies the properties in Theorem \ref{theorem1.6}, $T_m$ defines a bounded operator in $L^p(\mathbb R^n)$, for every $1<p<\infty$.

\begin{remark}\label{rem}
Note that, since the spectrum $\sigma(H)$ of $H$ is contained in $[1,\infty)$, Theorems \ref{nteo1} to \ref{theorem1.6} can be also established when the Hermite operator $H$ is replaced by the operator $H-\omega$, where $\omega <1$. In order to simplify the proofs we prefer to state the results for the Hermite operator, that is, in the case $\omega=0$.
\end{remark}

Bongioanni and Torrea (\cite{BT1} and \cite{BT2}) studied Sobolev spaces in the Hermite setting. In this paper we consider Hermite Sobolev spaces in the Banach valued context.

The Hermite operator $\tilde H$ admits the following factorization
$$
\tilde H = \frac 12 \sum_{j=1}^n (A_j A_{-j} + A_{-j} A_j) ,
$$
where $A_j = \frac{d}{dx_j} +x_j$ and $A_{-j} = -\frac{d}{dx_j} +x_j$, $j=1,\dots,n$.

Let $\mathbb{B}$ be a Banach space, $1<p<\infty$ and $\ell\in \mathbb{N}\setminus\{0\}$. The Hermite Sobolev $W_{H,\ell}^p (\mathbb{R}^n,\mathbb{B})$ is constituted by all the functions $f\in L^p (\mathbb{R}^n,\mathbb{B})$ such that, for every $j_i \in \mathbb{Z}$, $1\le |j_i|\le n$, $1\le i \le m\le \ell$, $A_{j_1} \cdots A_{j_m} f \in L^p (\mathbb{R}^n,\mathbb{B})$ in a distributional sense. The norm $\Vert \cdot \Vert_{W_{H,\ell}^p (\mathbb{R}^n,\mathbb{B})}$ in $W_{H,\ell}^p (\mathbb{R}^n,\mathbb{B})$ is defined by
$$
\Vert f \Vert_{W_{H,\ell}^p (\mathbb{R}^n,\mathbb{B})} = \Vert f \Vert_{L^p (\mathbb{R}^n,\mathbb{B})} +\sum_{\substack{j_i\in\mathbb{Z},\; 1\le |j_i|\le n\\ i=1,\dots,m\le \ell}} \Vert A_{j_1} \cdots A_{j_m} f \Vert_{L^p (\mathbb{R}^n,\mathbb{B})}, \;\; f\in W_{H,\ell}^p (\mathbb{R}^n,\mathbb{B}).
$$
We denote by $\tilde W_{H,\ell}^p (\mathbb{R}^n,\mathbb{B})$ the Hermite Sobolev space that consists of all those functions $f\in L^p (\mathbb{R}^n,\mathbb{B})$ such that, for every $j_i\in\mathbb{N}$, $1\le j_i\le n$, $1\le i \le m \le \ell$, $A_{-j_1} \cdots A_{-j_m} f \in L^p (\mathbb{R}^n,\mathbb{B})$ in a distributional sense. The norm $\Vert \cdot \Vert_{\tilde W_{H,\ell}^p (\mathbb{R}^n,\mathbb{B})}$ in $\tilde W_{H,\ell}^p (\mathbb{R}^n,\mathbb{B})$ is given by
$$
\Vert f \Vert_{\tilde W_{H,\ell}^p (\mathbb{R}^n,\mathbb{B})} = \Vert f \Vert_{L^p (\mathbb{R}^n,\mathbb{B})} +\sum_{\substack{j_i\in\mathbb{N},\; 1\le j_i\le n \\ i=1,\dots,m\le \ell}} \Vert A_{-j_1} \cdots A_{-j_m} f \Vert_{L^p (\mathbb{R}^n,\mathbb{B})}, \;\; f\in \tilde W_{H,\ell}^p (\mathbb{R}^n,\mathbb{B}).
$$
Let $\beta>0$. The $-\beta$-power $H^{-\beta}$ of the Hermite operator is defined by
$$
H^{-\beta} f= \frac{1}{\Gamma(\beta)}\int_0^\infty W_t^H(f)t^{\beta-1}dt,\;\; f\in L^p (\mathbb{R}^n,\mathbb{B}),\; 1<p<\infty.
$$
Since $H^{-\beta}$ is a positive operator, by \cite[Theorem 1]{BT1}, the operator $H^{-\beta}$ is bounded and one to one from $L^p (\mathbb{R}^n,\mathbb{B})$ into itself, for every $1<p<\infty$.

The potential space $L_{H,\beta}^p (\mathbb{R}^n,\mathbb{B})$ is constituted by all those $f\in L^p (\mathbb{R}^n,\mathbb{B})$ such that $f = H^{-\beta} g$, for some $g\in L^p (\mathbb{R}^n,\mathbb{B})$. The norm $\Vert \cdot \Vert_{L_{H,\beta}^p (\mathbb{R}^n,\mathbb{B})}$ on $L_{H,\beta}^p (\mathbb{R}^n,\mathbb{B})$ is given by
$$
\Vert f \Vert_{L_{H,\beta}^p (\mathbb{R}^n,\mathbb{B})} = \Vert g \Vert_{L^p (\mathbb{R}^n,\mathbb{B})},
$$
for every $f= H^{-\beta} g$, being $g\in L^p (\mathbb{R}^n,\mathbb{B})$.

The following result is a vector valued version of \cite[Theorem 4]{BT1} and \cite[Theorem 1]{BT2}, and Theorem \ref{theorem1.6} plays an important role in our proof.
\begin{teo}\label{theorem1.7}
Suppose that $\mathbb{B}$ is isomorphic to a closed subquotient of $[X,\mathcal{Q}]_\theta$ where $\theta \in (0,1)$, $X$ is a UMD Banach space, and $\mathcal{Q}$ is a Hilbert space, and that $\mathbb{B}$ has the property ($\alpha$). Then, for every $\ell\in \mathbb{N}\setminus\{0\}$,
$$
\tilde W_{H,\ell}^p (\mathbb{R}^n,\mathbb{B}) = W_{H,\ell}^p (\mathbb{R}^n,\mathbb{B}) = L_{H,\ell}^p (\mathbb{R}^n,\mathbb{B}),
$$
provided that $\left|\frac{2}{p}-1\right|<\theta$.
\end{teo}

As it was mentioned, Segovia and Wheeden (\cite{SW}) characterized classical Sobolev spaces by using square functions associated to Poisson semigroup involving fractional derivatives. In \cite{BFRTT} Hermite Sobolev spaces are described by fractional Littlewood-Paley functions defined by Hermite-Poisson semigroup. We now extend the results in \cite{BFRTT} for the Hermite Sobolev spaces to a Banach valued setting.

Let $\mathbb{B}$ be a Banach space and $1<p<\infty$. Suppose that $\beta >0$ and $k\in\mathbb{N}\setminus\{0\}$. We consider the Littlewood-Paley operator $G_{P,\beta,k;\mathbb{B}}^H$ defined by
$$
G_{P,\beta,k;\mathbb{B}}^H (f)(t,x)=t^\beta \partial_t^k P_t^H (f)(x,t),\;\;x\in\mathbb{R}^n~\textup{and}~t>0,
$$
for every $f\in L^p(\mathbb{R}^n,\mathbb{B})$.

By $F_{\beta,k}^H (\mathbb{R}^n,\mathbb{B})$ we denote the space that consists of all those functions $f\in L^p(\mathbb{R}^n,\mathbb{B})$ such that $G^H_{P,\beta,k;\mathbb{B}}(f) \in L^p(\mathbb{R}^n,\gamma(\mathcal{H}^n,\mathbb{B}))$. $F_{\beta,k}^H (\mathbb{R}^n,\mathbb{B})$ is endowed with the norm $\Vert \cdot \Vert_{F_{\beta,k}^H (\mathbb{R}^n,\mathbb{B})}$ defined by
$$
\Vert f \Vert_{F_{\beta,k}^H (\mathbb{R}^n,\mathbb{B})} = \Vert f \Vert_{L^p(\mathbb{R}^n,\mathbb{B})} + \Vert G^H_{P,\beta,k;\mathbb{B}}(f) \Vert_{L^p(\mathbb{R}^n,\gamma(\mathcal{H}^n,\mathbb{B}))},\;\; f \in F_{\beta,k}^H (\mathbb{R}^n,\mathbb{B}).
$$
Note that the space $F_{\beta,k}^H (\mathbb{R}^n,\mathbb{B})$ can be considered as a Triebel-Lizorkin type space in the Hermite setting.

\begin{teo}\label{nteo8}
Let $\mathbb{B}$ be a UMD Banach space and $1<p<\infty$. Suppose that $\beta >0$ and $k\in\mathbb{N}$ is such that $k>\beta$. Then, $L_{H,\beta}^p (\mathbb{R}^n,\mathbb{B}) = F_{k-\beta,k}^H (\mathbb{R}^n,\mathbb{B})$.
\end{teo}

In the following sections of this paper we present proofs for the Theorems. Throughout this paper by $C$ and $c$ we always represent positive constants but not necessarily the same in each occurrence.

The authors would like to thank Professor J.L. Torrea (UAM, Madrid). He has told us about Meda's multiplier theorem several years ago.

\section{Proof of Theorem \ref{theorem1.4}}

Let $1<p<\infty$. For every $\alpha=(\alpha_1,\dots,\alpha_n) \in (\mathbb{N}\setminus\{0\})^n$ we consider the operator $G_{P,\alpha;\mathbb{B}}^H$ defined by
\begin{equation}\label{n3.1}
G_{P,\alpha;\mathbb{B}}^H (f)(t,x) = \int_{\mathbb{R}^n} \prod_{j=1}^n t_j^{\alpha_j} \partial_{t_j}^{\alpha_j} P_{t_j}^H (x_j,y_j) f(y) dy,
\end{equation}
where $t=(t_1,\dots,t_n)\in (0,\infty)^n$ and $x=(x_1,\dots,x_n)\in\mathbb{R}^n$, for every $f\in L^p(\mathbb{R}^n,\mathbb{B})$.

In order to show that $G_{P,\alpha;\mathbb{B}}^H$ is a bounded operator from $L^p(\mathbb{R}^n,\mathbb{B})$ into $L^p(\mathbb{R}^n,\gamma(\mathcal{H}^n,\mathbb{B}))$ we proceed by induction on the dimension $n$.

According to \cite[Theorem 1]{BCCFR1} if $\alpha\in\mathbb{N}\setminus\{0\}$, the operator
$$
G_{P,\alpha;\mathbb{B}}^H (f)(t,x) = \int_{\mathbb{R}} t^\alpha \partial_t^\alpha P_t^H (x,y) f(y) dy,\;\; t\in (0,\infty)~\textup{and}~x\in\mathbb{R},
$$
is bounded from $L^p(\mathbb{R},\mathbb{B})$ into $L^p(\mathbb{R},\gamma(\mathcal{H}^1,\mathbb{B}))$

Suppose now that $n\in\mathbb{N}\setminus\{0\}$ and that for every $\alpha\in(\mathbb{N}\setminus\{0\})^n$ the operator defined by (\ref{n3.1}) is bounded from $L^p(\mathbb{R}^n,\mathbb{B})$ into $L^p(\mathbb{R}^n,\gamma(\mathcal{H}^n,\mathbb{B}))$. Let now $\beta =(\beta_1,\dots,\beta_{n+1})\in(\mathbb{N}\setminus\{0\})^{n+1}$. We are going to see that the operator defined through
$$
G_{P,\beta;\mathbb{B}}^H (f)(t,x) = \int_{\mathbb{R}^{n+1}} \prod_{j=1}^{n+1} t_j^{\beta_j} \partial_{t_j}^{\beta_j} P_{t_j}^H (x_j,y_j) f(y) dy,\;\; f\in L^p(\mathbb{R}^{n+1},\mathbb{B}),
$$
where $t=(t_1,\dots,t_{n+1}) \in (0,\infty)^{n+1}$ and $x=(x_1,\dots,x_{n+1})\in\mathbb{R}^{n+1}$, is bounded from $L^p(\mathbb{R}^{n+1},\mathbb{B})$ into $L^p(\mathbb{R}^{n+1},\gamma(\mathcal{H}^{n+1},\mathbb{B}))$.

Assume firstly that $f\in L^p(\mathbb{R})\otimes L^p(\mathbb{R}^n)\otimes \mathbb{B}$, that is $f=\sum_{i=1}^k g_i v_i b_i$, where $g_i \in L^p(\mathbb{R})$, $v_i\in L^p(\mathbb{R}^n)$ and $b_i \in\mathbb{B}$, $i=1,\dots,k \in \mathbb{N}\setminus\{0\}$. As it is well-known the space $L^p(\mathbb{R})\otimes L^p(\mathbb{R}^n)\otimes \mathbb{B}$ is dense in $L^p(\mathbb{R}^{n+1},\mathbb{B})$. It is clear that
$$
G_{P,\beta;\mathbb{B}}^H (f)(t,x) = \sum_{i=1}^k G_{P,\beta_1;\mathbb{C}}^H (g_i)(t_1,x_1) G_{P,\tilde\beta;\mathbb{C}}^H(v_i)(\tilde t,\tilde x) b_i,
$$
where $t=(t_1,t_2,\dots,t_{n+1})\in (0,\infty)^{n+1}$, $\tilde t=(t_2,\dots,t_{n+1})\in (0,\infty)^n$, $x=(x_1,x_2,\dots,x_{n+1})\in \mathbb{R}^{n+1}$, $\tilde x = (x_2,\dots,x_{n+1})\in\mathbb{R}^n$, and $\tilde \beta=(\beta_2,\dots,\beta_{n+1})\in(\mathbb{N}\setminus\{0\})^n$.

According to the induction hypothesis, if $v\in L^p(\mathbb{R}^n)$, for almost every $x\in\mathbb{R}^n$, $G_{P,\tilde\beta;\mathbb{C}}^H (v)(\cdot,x)\in \mathcal{H}^n$. Moreover, if $g\in L^p(\mathbb{R})$, for almost every $x\in\mathbb{R}$, $G_{P,\beta_1;\mathbb{C}}^H (g)(\cdot,x)\in \mathcal{H}^1$. Hence, we can choose a subset $A_1$ of $\mathbb{R}$ and a subset $A_2$ of $\mathbb{R}^n$ such that $|\mathbb{R}\setminus A_1|=0$, $|\mathbb{R}^n\setminus A_2|=0$, and, for every $i=1,\dots,k$ and $(x_1,\tilde x) \in A_1\times A_2$,
$$
G_{P,\beta_1;\mathbb{C}}^H (g_i)(\cdot,x_1) G_{P,\tilde \beta;\mathbb{C}}^H (v_i)(\cdot,\tilde x) \in \mathcal{H}^{n+1}.
$$
We conclude that, for every $x\in A_1\times A_2$, $G_{P,\beta;\mathbb{B}}^H (f)(\cdot,x) \in \mathcal{L}(\mathcal{H}^{n+1},\mathbb{B})$. Note that actually, for every $x\in A_1\times A_2$, $G_{P,\beta;\mathbb{B}}^H (f)(\cdot,x) \in \mathcal{L}(\mathcal{H}^{n+1},\mathbb{B}_k)$, where $\mathbb{B}_k=\langle\{b_i\}_{i=1}^k\rangle$ is the linear space generated by $\{b_i\}_{i=1}^k$. Since $\gamma(\mathcal{H}^l,\mathbb{C}) = \mathcal{H}^l$, $l\in \mathbb{N}\setminus\{0\}$, we get, for every $x\in A_1\times A_2$, $G_{P,\beta;\mathbb{B}}^H (f)(\cdot,x)\in \gamma(\mathcal{H}^{n+1},\mathbb{B}_k)$. We consider the function
$$
\begin{array}{cccl}
F: & \mathbb{R}^{n+1} & \rightarrow & \gamma(\mathcal{H}^{n+1},\mathbb{B})\\
        & x & \mapsto & F(x) = \begin{cases}
                                    G_{P,\beta;\mathbb{B}}^H (f)(\cdot,x),& x\in A_1\times A_2\\
                                    0,& x \in \mathbb{R}^{n+1} \setminus (A_1\times A_2).
                                \end{cases}
\end{array}
$$
According to \cite[Lemma 2.5]{vNVW} $F$ is strongly measurable if, and only if, for every $h\in \mathcal{H}^{n+1}$ the function
$$
\begin{array}{cccl}
F_h: & \mathbb{R}^{n+1} & \rightarrow & \mathbb{B}\\
        & x & \mapsto & [F(x)](h) = \begin{cases}
                                        \int_{(0,\infty)^{n+1}} G_{P,\beta;\mathbb{B}}^H (f)(t,x) h(t) \frac{dt_1\cdots dt_{n+1}}{t_1\cdots t_{n+1}},& x\in A_1\times A_2\\
                                        0,& x \in \mathbb{R}^{n+1} \setminus (A_1\times A_2)
                                    \end{cases}
\end{array}
$$
is strongly measurable. Let $h\in \mathcal{H}^{n+1}$. We have that
\begin{multline*}
\int_{(0,\infty)^{n+1}} G_{P,\beta;\mathbb{B}}^H (f)(t,x) h(t) \frac{dt_1\cdots dt_{n+1}}{t_1\cdots t_{n+1}} =\\ \sum_{i=1}^k \int_{(0,\infty)^{n+1}} G_{P,\beta_1;\mathbb{C}}^H (g_i)(t_1,x_1) G_{P,\tilde\beta;\mathbb{C}}^H (v_i) (\tilde t, \tilde x) h(t) \frac{dt_1\cdots dt_{n+1}}{t_1\cdots t_{n+1}} b_i,
\end{multline*}
where $x=(x_1,\tilde x)\in A_1\times A_2$. It is clear that the range of $F_h$ is contained in $\mathbb{B}_k$. Hence, the range of $F_h$ is separable subset of $\mathbb{B}$.

Assume now that $S\in\mathbb{B}^*$. We can write
$$
\langle S, F_h(x)\rangle_{\mathbb{B}^*,\mathbb{B}} =
\begin{cases}
\sum_{i=1}^k \int_{(0,\infty)^{n+1}} G_{P,\beta;\mathbb{C}}^H (f_i)(t,x) h(t) \frac{dt_1\cdots dt_{n+1}}{t_1\cdots t_{n+1}} \langle S, b_i\rangle_{\mathbb{B}^*,\mathbb{B}},& x\in A_1\times A_2\\
0,& x\in \mathbb{R}^{n+1}\setminus(A_1\times A_2),
\end{cases}
$$
where $f_i(y)=g_i(y_1)v_i(\tilde y)$, $y=(y_1,\dots,y_{n+1})\in \mathbb{R}^{n+1}$, $\tilde y =(y_1,\dots,y_{n+1})\in \mathbb{R}^n$.

The function $\langle S, F_h(x)\rangle_{\mathbb{B}^*,\mathbb{B}}$, $x\in\mathbb{R}^{n+1}$, is Lebesgue measurable in $\mathbb{R}^{n+1}$, because, for $i=1,\dots,k$, the function
$$
\begin{array}{ccl}
 \mathbb{R}^{n+1} & \rightarrow & \mathcal{H}^{n+1}\\
         x & \mapsto & G_{P,\beta;\mathbb{C}}^H (f_i)(x)
\end{array}
$$
is strongly measurable (this property is implicitly included in \cite[Theorem 2.4]{Wr}).

We conclude that the function $F$ is strongly measurable.

Since $\mathbb{B}$ has the property ($\alpha$) we have that $\gamma(\mathcal{H}^{n+1},\mathbb{B}) \simeq \gamma(\mathcal{H}^1,\gamma(\mathcal{H}^n,\mathbb{B}))$ (\cite[Corollary 3.5]{vNW}). By using the induction hypothesis and \cite[Theorem 1]{BCCFR1} and taking into account that $\gamma(\mathcal{H}^l,\mathbb{B})$, $l\in\mathbb{N}\setminus\{0\}$, is UMD with the property ($\alpha$) we obtain
\begin{eqnarray*}
\lefteqn{\Vert G_{P,\beta;\mathbb{B}}^H (f) \Vert_{L^p(\mathbb{R}^{n+1},\gamma(\mathcal{H}^{n+1},\mathbb{B}))}}\\
&=& \left(\int_{\mathbb{R}^{n+1}} \Vert G_{P,\beta;\mathbb{B}}^H (f) (\cdot,x) \Vert_{\gamma(\mathcal{H}^{n+1},\mathbb{B})}^p dx\right)^{\frac{1}{p}}\\
&=& \left( \int_{\mathbb{R}^{n+1}} \left\Vert \int_{\mathbb{R}} t_1^{\beta_1} \partial_{t_1}^{\beta_1} P_{t_1}^H (x_1,y_1) G_{P,\tilde \beta;\mathbb{B}}^H (f)(\tilde t,\tilde x) dy_1\right\Vert^p_{\gamma(\mathcal{H}^1,\gamma(\mathcal{H}^n,\mathbb{B}))} dx \right)^{\frac{1}{p}}\\
&\le& C \left( \int_{\mathbb{R}} \int_{\mathbb{R}^n} \Vert G_{P,\tilde \beta;\mathbb{B}}^H (f(x_1,\tilde y))(\tilde t,\tilde x) \Vert_{\gamma(\mathcal{H}^n,\mathbb{B})}^p d\tilde x dx_1 \right)^{\frac{1}{p}}\\
&\le& C \Vert f \Vert_{L^p(\mathbb{R}^{n+1})}
\end{eqnarray*}
Since $L^p(\mathbb{R})\otimes L^p(\mathbb{R}^n)\otimes \mathbb{B}$ is a dense subspace of $L^p(\mathbb{R}^{n+1},\mathbb{B})$, the operator $G_{P,\beta;\mathbb{B}}^H$ can be extended to $L^p(\mathbb{R}^{n+1},\mathbb{B})$ as a bounded operator from $L^p(\mathbb{R}^{n+1},\mathbb{B})$ into $L^p(\mathbb{R}^{n+1},\gamma(\mathcal{H}^{n+1},\mathbb{B}))$. This extension operator is denoted by $\tilde{G}_{P,\beta;\mathbb{B}}^H$.

We are going to show that, for every $f\in L^p(\mathbb{R}^{n+1},\mathbb{B})$,
$$
\tilde{G}_{P,\beta;\mathbb{B}}^H (f)(x) = G_{P,\beta;\mathbb{B}}^H (f)(\cdot, x),\;\;\textup{a.e.}~x\in\mathbb{R}^{n+1}.
$$
Let $f\in L^p(\mathbb{R}^{n+1},\mathbb{B})$. We choose a sequence $\{f_m\}_{m=1}^\infty \subset S(\mathbb{R}) \otimes S(\mathbb{R}^n) \otimes \mathbb{B}$ such that $f_m \to f$, as $m\to\infty$, in $L^p(\mathbb{R}^{n+1},\mathbb{B})$. Then, $G_{P,\beta;\mathbb{B}}^H (f_m) \to \tilde{G}_{P,\beta;\mathbb{B}}^H (f)$, as $m\to\infty$ in $L^p(\mathbb{R}^{n+1},\gamma(\mathcal{H}^{n+1},\mathbb{B}))$, and there exists a subsequence $\{f_{m_k}\}_{k=1}^\infty$ of $\{f_m\}_{m=1}^\infty$ such that $G_{P,\beta;\mathbb{B}}^H (f_{m_k})(\cdot,x) \to \tilde{G}_{P,\beta;\mathbb{B}}^H (f)(x)$, as $k\to \infty$ in $\gamma(\mathcal{H}^{n+1},\mathbb{B})$, for every $x\in A\subset \mathbb{R}^{n+1}$ being $|\mathbb{R}^{n+1}\setminus A|=0$. Since $\gamma(\mathcal{H}^{n+1},\mathbb{B})$ is continuously contained in the space $\mathcal{L}(\mathcal{H}^{n+1},\mathbb{B})$, we have that $G_{P,\beta;\mathbb{B}}^H (f_m)(\cdot,x) \to \tilde{G}_{P,\beta;\mathbb{B}}^H (f)(x)$, as $k\to \infty$, in $\mathcal{L}(\mathcal{H}^{n+1},\mathbb{B})$, for every $x\in A$.

On the other hand, for every $g\in L^p(\mathbb{R}^{n+1},\mathbb{B})$ and  $h\in L^2\left((\delta,\infty)^{n+1},\frac{dt_1\cdots dt_{n+1}}{t_1\cdots t_{n+1}}\right)$, where $\delta>0$, we get
\begin{equation}\label{na0}
\begin{split}
\lefteqn{\left\Vert \int_{(0,\infty)^{n+1}} G_{P,\beta;\mathbb{B}}^H (g)(t,x) h(t) \frac{dt_1\cdots dt_{n+1}}{t_1\cdots t_{n+1}}\right\Vert_{\mathbb{B}}}\\
& \le \left( \int_{(\delta,\infty)^{n+1}} \Vert G_{P,\beta;\mathbb{B}}^H (g)(t,x) \Vert_{\mathbb{B}}^2 \frac{dt_1\cdots dt_{n+1}}{t_1\cdots t_{n+1}} \right)^{\frac12} \Vert h\Vert_{\mathcal{H}^{n+1}}\\
& \le \int_{\mathbb{R}^{n+1}} \Vert g(y)\Vert_{\mathbb{B}} \left ( \prod_{j=1}^{n+1} \int_\delta^\infty \left| t_j^{\beta_j} \partial_{t_j}^{\beta_j} P_{t_j}^H (x_j,y_j) \right|^2 \frac{dt_j}{t_j}\right)^{\frac12} dy \Vert h\Vert_{\mathcal{H}^{n+1}}.
\end{split}
\end{equation}

Let $l\in\mathbb{N}$, $l\ge 1$. The subordination formula allows us to write,
\begin{equation}\label{na1}
\begin{split}
    \partial_s^l P^H_s(v,z) &= \partial_s^l \left[ \frac{s}{\sqrt{4\pi}} \int_0^\infty \frac{e^{-\frac{s^2}{4u}}}{u^\frac32} W^H_u(v,z) du \right]\\
    &= -\frac{1}{\sqrt{\pi}} \int_0^\infty \partial_s^{l +1} \left( \frac{e^{-\frac{s^2}{4u}}}{u^\frac12} \right) W^H_u(v,z) du, \;\; z, v \in \mathbb{R}~\textup{and}~s\in(0,\infty).
\end{split}
\end{equation}
According to Fa\`{a} di Bruno's formula (\cite[(4.6)]{BGLNU}) we have that,
\begin{equation}\label{na2}
\begin{split}
\partial_s^{l +1} \left( e^{-\frac{s^2}{4u}} \right) &= \frac{1}{2^{l+1}u^{\frac{l+1}{2}}} \partial_v^{l +1} \left( e^{-v^2} \right)_{\vert v= \frac{s}{2\sqrt{u}}}\\
&= \frac{1}{2^{l+1}u^{\frac{l+1}{2}}} \sum_{k=0}^{\frac{l+1}{2}} (-1)^{l - k} E_{l+1,k} \left( \frac{s}{2\sqrt{u}} \right)^{l+1-2k} e^{-\frac{s^2}{4u}},\;\; s,u \in (0,\infty),
\end{split}
\end{equation}
where $E_{l+1,k} = \frac{2^{l+1-2k} (l+1)!}{k! (l+1-2k)!}$, $0\le k \le \frac{l+1}{2}$.

From (\ref{na2}) we deduce that
$$
    \left\vert \partial_s^{l +1} \left( e^{-\frac{s^2}{4u}} \right) \right\vert \le C \frac{e^{-\frac{cs^2}{u}}}{u^{\frac{l+1}{2}}},\;\; s,u \in (0,\infty).
$$
The equality (\ref{na1}), since $W_u^H(v,z) \le  \frac{e^{-\frac{|v-z|^2}{4u}}}{\sqrt{2\pi u}}$, $v\in\mathbb{R}$ and $u\in(0,\infty)$ \cite[(2.9)]{StemTo1}, leads to
$$
|\partial_s^l P^H_s(v,z)| \le \int_0^\infty \frac{e^{-\frac{-c(s^2+|v-z|^2)}{u}}}{u^{\frac{l+3}{2}}} du \le C \frac{1}{(s+|v-z|)^{1+l}},\;\; v,z \in \mathbb{R}~\textup{and}~s\in(0,\infty).
$$
Hence we get
\begin{align*}
\int_\delta^\infty \left| s^l \partial_s^l P_s^H(v,z) \right|^2 \frac{ds}{s} &\le C \int_\delta^\infty \frac{s^{2l-1}}{(s+|v-z|)^{2l+2}} ds\\
&\le C \int_\delta^\infty \frac{ds}{(s+|v-z|)^3} = \frac{C}{(|v-z|+\delta)^2},\;\; v,z \in \mathbb{R}.
\end{align*}
From (\ref{na0}) it follows that, for every $h\in L^2\left((\delta,\infty)^{n+1},\frac{dt_1\cdots dt_{n+1}}{t_1\cdots t_{n+1}}\right)$,
\begin{equation}\label{na3}
\begin{split}
\lefteqn{\left\Vert \int_{(0,\infty)^{n+1}} G_{P,\beta;\mathbb{B}}^H (g)(t,x) h(t) \frac{dt_1\cdots dt_{n+1}}{t_1\cdots t_{n+1}} \right\Vert_{\mathbb{B}}}\hspace{2cm}\\
&\le \int_{\mathbb{R}^{n+1}} \Vert g(y) \Vert_{\mathbb{B}} \prod_{j=1}^{n+1} \frac{1}{(\delta + |x_j-y_j|)} dy \Vert h \Vert_{\mathcal{H}^{n+1}}\\
&\le C \Vert g \Vert_{L^p(\mathbb{R}^{n+1},\mathbb{B})} \Vert h \Vert_{\mathcal{H}^{n+1}}.
\end{split}
\end{equation}
Let $h\in \mathcal{H}^{n+1}$ such that $\textup{supp}~h \subset (\delta,\infty)^{n+1}$, for a certain $\delta>0$. From (\ref{na3}) we deduce that
$$
\int_{(0,\infty)^{n+1}} G_{P,\beta;\mathbb{B}}^H (f_{m_k})(t,x) h(t) \frac{dt_1\cdots dt_{n+1}}{t_1\cdots t_{n+1}} \to \int_{(0,\infty)^{n+1}} G_{P,\beta;\mathbb{B}}^H (f)(t,x) h(t) \frac{dt_1\cdots dt_{n+1}}{t_1\cdots t_{n+1}}
$$
in $\mathbb{B}$, for every $x\in\mathbb{R}^{n+1}$, and, for every $S\in\mathbb{B}^*$,
\begin{eqnarray*}
\lefteqn{\left\langle S, \int_{(0,\infty)^{n+1}} G_{P,\beta;\mathbb{B}}^H (f_{m_k})(t,x) h(t) \frac{dt_1\cdots dt_{n+1}}{t_1\cdots t_{n+1}} \right\rangle_{\mathbb{B}^*,\mathbb{B}}}\\
&=&\int_{(0,\infty)^{n+1}} \langle S, G_{P,\beta;\mathbb{B}}^H (f_{m_k})(t,x)\rangle_{\mathbb{B}^*,\mathbb{B}} h(t) \frac{dt_1\cdots dt_{n+1}}{t_1\cdots t_{n+1}}\\
&\to& \int_{(0,\infty)^{n+1}} \langle S, G_{P,\beta;\mathbb{B}}^H (f)(t,x)\rangle_{\mathbb{B}^*,\mathbb{B}} h(t) \frac{dt_1\cdots dt_{n+1}}{t_1\cdots t_{n+1}} \\
&=& \left\langle S, \int_{(0,\infty)^{n+1}} G_{P,\beta;\mathbb{B}}^H (f)(t,x) h(t) \frac{dt_1\cdots dt_{n+1}}{t_1\cdots t_{n+1}} \right\rangle_{\mathbb{B}^*,\mathbb{B}}
\end{eqnarray*}
as $k\to\infty$, for every $x\in \mathbb{R}^{n+1}$.

Also, for every $S\in \mathbb{B}^*$, and $x\in A$, we have that
\begin{multline*}
\langle S, G_{P,\beta;\mathbb{B}}^H (f_{m_k})(\cdot,x)[h]\rangle_{\mathbb{B}^*,\mathbb{B}} = \int_{(0,\infty)^{n+1}} \langle S, G_{P,\beta;\mathbb{B}}^H (f_{m_k})(t,x)\rangle_{\mathbb{B}^*,\mathbb{B}} h(t) \frac{dt_1\cdots dt_{n+1}}{t_1\cdots t_{n+1}}\\ \to \langle S, \tilde{G}_{P,\beta;\mathbb{B}}^H (f)(x)[h]\rangle_{\mathbb{B}^*,\mathbb{B}},
\end{multline*}
as $k\to\infty$. Hence, for every $S\in \mathbb{B}^*$,
$$
\langle S, \tilde{G}_{P,\beta;\mathbb{B}}^H (f)(x)[h]\rangle_{\mathbb{B}^*,\mathbb{B}}= \int_{(0,\infty)^{n+1}} \langle S, G_{P,\beta;\mathbb{B}}^H (f)(t,x)\rangle_{\mathbb{B}^*,\mathbb{B}} h(t) \frac{dt_1\cdots dt_{n+1}}{t_1\cdots t_{n+1}},\;\; x\in A,
$$
and
$$
\left| \int_{(0,\infty)^{n+1}} \langle S, G_{P,\beta;\mathbb{B}}^H (f)(t,x)\rangle_{\mathbb{B}^*,\mathbb{B}} h(t) \frac{dt_1\cdots dt_{n+1}}{t_1\cdots t_{n+1}} \right|
\le C ||S||_{\mathbb{B}^*} ||G_{P,\beta;\mathbb{B}}^H (f)(x)||_{\mathcal{L}(\mathcal{H}^{n+1},\mathbb{B})} ||h||_{\mathcal{H}^{n+1}} ,
$$
We conclude that, for every $S\in \mathbb{B}^*$ and $x\in A$, $\langle S,G_{P,\beta;\mathbb{B}}^H(f)(\cdot,x) \rangle \in \mathcal{H}^{n+1}$ and
$$
\int_{(0,\infty)^{n+1}} \langle S, G_{P,\beta;\mathbb{B}}^H (f)(t,x)\rangle_{\mathbb{B}^*,\mathbb{B}} h(t) \frac{dt_1\cdots dt_{n+1}}{t_1\cdots t_{n+1}} =
\langle S, \tilde{G}_{P,\beta;\mathbb{B}}^H (f)(x)[h]\rangle, \;\; h\in \mathcal{H}^{n+1}.
$$
Then,
$$
G_{P,\beta;\mathbb{B}}^H (f)(\cdot,x) = \tilde{G}_{P,\beta;\mathbb{B}}^H (f)(x),\;\; x\in A,
$$
as elements of $\mathcal{L}(\mathcal{H}^{n+1},\mathbb{B})$

Thus, the induction is completed and we prove that the operator $G_{P,\beta;\mathbb{B}}^H$ is bounded from $L^p(\mathbb{R}^{n+1},\mathbb{B})$ into $L^p(\mathbb{R}^{n+1},\gamma(\mathcal{H}^{n+1},\mathbb{B}))$, for every $\beta \in (\mathbb{N}\setminus\{0\})^{n+1}$ and $n\in\mathbb{N}$.

Our next objective is to see that there exists $C>0$ such that
$$
\Vert f \Vert_{L^p(\mathbb{R}^{n+1},\mathbb{B})} \le C \Vert G_{P,\alpha;\mathbb{B}}^H \Vert_{L^p(\mathbb{R}^{n+1},\gamma(\mathcal{H}^{n+1},\mathbb{B}))}, \;\; f\in L^p(\mathbb{R}^{n+1},\mathbb{B}).
$$
By using standard spectral arguments we can show that (see \cite[Proposition 2.1]{BCFR}), for every $f\in L^p(\mathbb{R}^{n+1})\otimes B$ and $g\in L^{p'}(\mathbb{R}^{n+1})\otimes B^*$, where $p'$ is the conjugated of $p$, that is $p'=\frac{p}{p-1}$,
\begin{equation}\label{pola}
\begin{split}
\lefteqn{\int_{\mathbb{R}^{n+1}} \int_{(0,\infty)^{n+1}} \langle G_{P,\alpha;\mathbb{B}^*}^H (g)(t,x),G_{P,\alpha;\mathbb{B}}^H (f)(t,x) \rangle_{\mathbb{B}^*,\mathbb{B}} \frac{dt_1\cdots dt_{n+1}}{t_1 \cdots t_{n+1}} dx} \hspace{4cm}\\
&= \prod_{j=1}^{n+1} \frac{\Gamma(2\alpha_j)}{2^{2\alpha_j}} \int_{\mathbb{R}^{n+1}} \langle g(x),f(x)\rangle_{\mathbb{B}^*,\mathbb{B}} dx
\end{split}
\end{equation}

Let $f\in L^p(\mathbb{R}^{n+1})\otimes B$. By \cite[Lemma 2.3]{GrafYang} and \cite[Proposition 2.2]{HW}, we have that
\begin{eqnarray*}
\lefteqn{||f||_{L^p(\mathbb{R}^{n+1},\mathbb{B})} = \sup_{\substack{g\in L^{p'}(\mathbb{R}^{n+1})\otimes \mathbb{B}^*\\ ||g||_{L^{p'}(\mathbb{R}^{n+1},\mathbb{B}^*)} \le 1}} \left| \int_{\mathbb{R}^{n+1}} \langle g(x), f(y) \rangle_{\mathbb{B}^*,\mathbb{B}} dy \right|}\\
&=& \sup_{\substack{g\in L^{p'}(\mathbb{R}^{n+1})\otimes \mathbb{B}^*\\ ||g||_{L^{p'}(\mathbb{R}^{n+1},\mathbb{B}^*)} \le 1}} \left| \prod_{j=1}^{n+1} \frac{2^{2\alpha_j}}{\Gamma(2\alpha_j)} \int_{\mathbb{R}^{n+1}} \int_{(0,\infty)^{n+1}} \langle G_{P,\alpha;\mathbb{B}^*}^H (g)(t,x), G_{P,\alpha;\mathbb{B}}^H (f)(t,x)\rangle_{\mathbb{B}^*,\mathbb{B}} \frac{dt_1\cdots dt_{n+1}}{t_1 \cdots t_{n+1}} dx \right|\\
&\le& \sup_{\substack{g\in L^{p'}(\mathbb{R}^{n+1})\otimes \mathbb{B}^*\\ ||g||_{L^{p'}(\mathbb{R}^{n+1},\mathbb{B}^*)} \le 1}} \int_{\mathbb{R}^{n+1}}  \Vert G_{P,\alpha;\mathbb{B}^*}^H (g)(\cdot,x) \Vert_{\gamma(\mathcal{H}^{n+1},\mathbb{B}^*)} \Vert G_{P,\alpha;\mathbb{B}}^H (f)(\cdot,x) \Vert_{\gamma(\mathcal{H}^{n+1},\mathbb{B})} dx.
\end{eqnarray*}
Since $\mathbb{B}^*$ is UMD and it has the property ($\alpha$), $G_{P,\alpha;\mathbb{B}^*}^H$ is bounded from $L^{p'}(\mathbb{R}^{n+1},\mathbb{B}^*)$ into $L^{p'}(\mathbb{R}^{n+1},\gamma(\mathcal{H}^{n+1},\mathbb{B}^*))$, and by using H\"older's inequality we get
$$
\Vert f \Vert_{L^p(\mathbb{R}^{n+1},\mathbb{B})} \le C \Vert G_{P,\alpha; \mathbb{B}}^H (f) \Vert_{L^p(\mathbb{R}^{n+1},\gamma(\mathcal{H}^{n+1},\mathbb{B}))}
$$
Since $L^p(\mathbb{R}^{n+1})\otimes \mathbb{B}$ is a dense subspace of $L^p(\mathbb{R}^{n+1},\mathbb{B})$ and $G_{P,\alpha;\mathbb{B}}^H$ is a bounded operator from $L^p(\mathbb{R}^{n+1},\mathbb{B})$ into $L^p(\mathbb{R}^{n+1},\gamma(\mathcal{H}^{n+1},\mathbb{B}))$ we conclude that
$$
\Vert f \Vert_{L^p(\mathbb{R}^{n+1},\mathbb{B})} \le C \Vert G_{P,\alpha; \mathbb{B}}^H (f) \Vert_{L^p(\mathbb{R}^{n+1},\gamma(\mathcal{H}^{n+1},\mathbb{B}))},\;\; f\in L^p(\mathbb{R}^{n+1},\mathbb{B}).
$$

\section{Proof of Theorem \ref{theorem1.5}}

Since $\textrm{m}\in L^\infty (\mathbb{R}^n)$ the multiplier operator $T_\textrm{m}$ defined by
$$
T_\textrm{m} (f) = \sum_{k_1,\dots,k_n=0}^\infty \textrm{m}(\lambda_{k_1},\dots,\lambda_{k_n}) c_{k_1,\dots,k_n} (f) h_{k_1,\dots,k_n}, \;\; f\in L^2(\mathbb{R}^n),
$$
is bounded from $L^2(\mathbb{R}^n)$ into itself. Here, for every $k_1,\dots,k_n \in \mathbb{N}$,
$$
h_{k_1,\dots,k_n} (x) = \prod_{j=1}^n h_{k_j} (x_j),\;\; x=(x_1,\dots,x_n)\in \mathbb{R}^n,
$$
and, for every $f\in L^2(\mathbb{R}^n)$,
$$
c_{k_1,\dots,k_n} (f) = \int_{\mathbb{R}^n} h_{k_1,\dots,k_n} (x) f(x) dx.
$$

Inspired by \cite{Me}, \cite{Wr} and \cite{Wr1} we consider, for every $\alpha=(\alpha_1,\dots,\alpha_n)\in\mathbb{N}^n$,
$$
m_\alpha (t,\lambda) = \prod_{j=1}^n (t_j \lambda_j)^{\alpha_j} e^{-\frac{t_j \lambda_j}{2}} M(\lambda),
$$
where $t=(t_1,\dots,t_n)\in (0,\infty)^n$, $\lambda=(\lambda_1,\dots,\lambda_n) \in \mathbb{R}^n$ and $M(\lambda)=\textrm{m}(\lambda_1^2,\dots,\lambda_n^2)$, $\lambda=(\lambda_1,\dots,\lambda_n)\in\mathbb{R}^n$, and
$$
\mathcal{M}_\alpha (t,u) = \int_{(0,\infty)^n} \prod_{j=1}^n \lambda_j^{-i u_j -1} m_\alpha (t,\lambda) d\lambda,
$$
with $t=(t_1,\dots,t_n)\in (0,\infty)^n$ and $u=(u_1,\dots,u_n)\in\mathbb{R}^n$.

It is clear that $T_\textrm{m} = M (\sqrt{H},\dots,\sqrt{H})$. By proceeding as in the proof of \cite[Lemma 3.2]{Wr} we can show that, for every $\alpha=(\alpha_1,\dots,\alpha_n)\in\mathbb{N}^n$ and $f\in L^2(\mathbb{R}^n)$,
\begin{equation}\label{4.1}
    G_{P,\alpha+1,\mathbb{B}}^H (T_\textrm{m} f)(t,x) = \frac{1}{2^n\pi^n} \int_{\mathbb{R}^n} \mathcal{M}_\alpha (t,u) \int_{\mathbb{R}^n} \prod_{j=1}^n t_j \partial_{t_j} P_{\frac{t_j}{2}}^H (x_j,y_j) L^{i \frac{u}{2}} f(y) dy du,
\end{equation}
with $t=(t_1,\dots,t_n)\in (0,\infty)^n$ and $x=(x_1,\dots,x_n)\in\mathbb{R}^n$ and where $\alpha+1=(\alpha_1+1,\dots,\alpha_n+1)$ and
$$
L^{i v} f = \sum_{k_1,\dots,k_n=0}^\infty \prod_{j=1}^n \lambda_{k_j}^{iv_j} c_{k_1,\dots,k_n} (f) h_{k_1,\dots,k_n}, \;\; f\in L^2(\mathbb{R}^n)~\textup{and}~v=(v_1,\dots,v_n).
$$
The imaginary power $H^{i\gamma}$, $\gamma\in \mathbb{R}$, of the Hermite operator $H$ is bounded from $L^p(\mathbb{R})$ into itself (\cite{BCCR} and \cite{StemTo3}). Thus, the operator $L^{iv}$ is bounded from $L^p(\mathbb{R}^n)$ into $L^p(\mathbb{R}^n)$ for every $v \in \mathbb{R}^n$.

Suppose that $f\in L^p(\mathbb{R}^n) \cap L^2(\mathbb{R}^n)$. Then, $T_\textrm{m} f \in L^2(\mathbb{R}^n)$. According to (\ref{4.1}), Minkowski's inequality leads to
$$
\Vert G_{P,\alpha+1;\mathbb{B}}^H (T_\textrm{m} f)(\cdot,x) \Vert_{\mathcal{H}^n} \le C \int_{\mathbb{R}^n} \sup_{t\in (0,\infty)^n} |\mathcal{M}_\alpha (t,u)| \Vert G_{P,1;\mathbb{B}}^H (L^{i \frac{u}{2}} f) (\cdot,x)\Vert_{\mathcal{H}^n} du
$$
Hence, since $G_{P,1;\mathbb{B}}^H$ is bounded from $L^p(\mathbb{R}^n)$ into $L^p(\mathbb{R}^n,\mathcal{H}^n)$
$$
\Vert G_{P,\alpha+1;\mathbb{B}}^H (T_\textrm{m} f)(\cdot,x) \Vert_{L^p(\mathbb{R}^n,\mathcal{H}^n)} \le C \int_{\mathbb{R}^n} \sup_{t\in (0,\infty)^n} |\mathcal{M}_\alpha (t,u)| \Vert L^{i \frac{u}{2}} \Vert_{L^p(\mathbb{R}^n) \to L^p(\mathbb{R}^n)} du \Vert f \Vert_{L^p(\mathbb{R}^n)}.
$$
By using the polarization formula (\ref{pola}) in the scalar case we get, for every $g\in L^{p'}(\mathbb{R}^n) \cap L^2(\mathbb{R}^n)$,
\begin{align*}
  \left|\int_{\mathbb{R}^n} T_\textrm{m}(f)(x) g(x) dx \right| &\le C \Vert G^H_{P,\alpha +1;\mathbb{B}}(T_\textrm{m} f)\Vert_{L^p(\mathbb{R}^n,\mathcal{H}^n)} \Vert G^H_{P,\alpha +1;\mathbb{B}}(g) \Vert_{L^{p'}(\mathbb{R}^n,\mathcal{H}^n)}\\
   &\le  C \Vert G^H_{P,\alpha +1;\mathbb{B}}(T_\textrm{m} f)\Vert_{L^p(\mathbb{R}^n,\mathcal{H}^n)} \Vert g \Vert_{L^{p'}(\mathbb{R}^n)}.
\end{align*}
Then, $T_\textrm{m}(f) \in L^p(\mathbb{R}^n)$ and
\begin{align*}
    \Vert T_\textrm{m}(f) \Vert_{L^p(\mathbb{R}^n)} &\le C \Vert G^H_{P,\alpha +1;\mathbb{B}}(T_\textrm{m} f)\Vert_{L^p(\mathbb{R}^n,\mathcal{H}^n)}\\
    & \le C \int_{\mathbb{R}^n} \sup_{t\in (0,\infty)^n} |\mathcal{M}_\alpha (t,u)| \Vert L^{i \frac{u}{2}} \Vert_{L^p(\mathbb{R}^n) \to L^p(\mathbb{R}^n)} du\Vert f \Vert_{L^p(\mathbb{R}^n)}.
\end{align*}
Hence, $T_\textrm{m}$ can be extended from $L^p(\mathbb{R}^n) \cap L^2(\mathbb{R}^n)$ to $L^p(\mathbb{R}^n)$ as a bounded operator from $L^p(\mathbb{R}^n)$ into itself.

Let now $f\in S(\mathbb{R}^n)\otimes \mathbb{B}$. Then, by defining $T_\textrm{m}$ on $L^p(\mathbb{R}^n)\otimes \mathbb{B}$ in the natural way, we have that $T_\textrm{m}(f) \in L^p(\mathbb{R}^n)\otimes \mathbb{B}$ and (\ref{4.1}) holds for $f$.

According to (\ref{1.6}) $\mathcal{M}_\alpha (\cdot,u) \in L^\infty ((0,\infty)^n)$, for every $u \in \Omega \subset \mathbb{R}^n$ such that $|\mathbb{R}^n\setminus \Omega|=0$. Hence, for every $u \in \Omega$, $\mathcal{M}_\alpha (\cdot,u)$ defines a pointwise multiplier in $\mathcal{H}^n$. By using the ideal property of $\gamma(\mathcal{H}^n,\mathbb{B})$ (\cite[Theorem 6.2]{VanNee}) we deduce that, for every $u\in \Omega$
\begin{align*}
\Vert \mathcal{M}_\alpha (\cdot,u) G_{P,1;\mathbb{B}}^H (L^{i \frac{u}{2}}f)(\cdot,x)\Vert_{\gamma(\mathcal{H}^n,\mathbb{B})} &\le C \Vert \mathcal{M}_\alpha (\cdot,u) \Vert_{\mathcal{H}^n \to \mathcal{H}^n} \Vert G_{P,1;\mathbb{B}}^H (L^{i \frac{u}{2}}f)(\cdot,x) \Vert_{\gamma(\mathcal{H}^n,\mathbb{B})}\\
&\le C \sup_{t\in (0,\infty)^n} |\mathcal{M}_\alpha (\cdot,u)| \Vert G_{P,1;\mathbb{B}}^H (L^{i \frac{u}{2}}f)(\cdot,x) \Vert_{\gamma(\mathcal{H}^n,\mathbb{B})}.
\end{align*}
Since $\mathbb{B}$ is a UMD Banach space the imaginary power $H^{i\gamma}$, $\gamma\in\mathbb{R}$, is bounded from $L^p(\mathbb{R},\mathbb{B})$ into itself (\cite{BCCR}). Then, the operator $L^{iv}$ is bounded from $L^p(\mathbb{R}^n,\mathbb{B})$ into itself, for every $v\in\mathbb{R}^n$. According to Theorem \ref{theorem1.4} we get, for every $u\in\Omega$,
\begin{eqnarray*}
\lefteqn{\Vert \mathcal{M}_\alpha (\cdot,u) G_{P,1;\mathbb{B}}^H (L^{i \frac{u}{2}}f) \Vert_{L^p(\mathbb{R}^n,\gamma(\mathcal{H}^n,\mathbb{B}))}}\\
&\le& C \sup_{t\in (0,\infty)^n} |\mathcal{M}_\alpha (\cdot,u)| \Vert G_{P,1;\mathbb{B}}^H (L^{i \frac{u}{2}}f)\Vert_{L^p(\mathbb{R}^n,\gamma(\mathcal{H}^n,\mathbb{B}))}\\
&\le& C \sup_{t\in (0,\infty)^n} |\mathcal{M}_\alpha (\cdot,u)| \Vert  L^{i \frac{u}{2}}f \Vert_{L^p(\mathbb{R}^n,\mathbb{B})}\\
&\le& C \sup_{t\in (0,\infty)^n} |\mathcal{M}_\alpha (\cdot,u)| \Vert  L^{i \frac{u}{2}} \Vert_{L^p(\mathbb{R}^n,\mathbb{B}) \to L^p(\mathbb{R}^n,\mathbb{B})} \Vert f \Vert_{L^p(\mathbb{R}^n,\mathbb{B})}.
\end{eqnarray*}
By using again Theorem \ref{theorem1.4} we obtain
\begin{align*}
    \Vert T_\textrm{m} f \Vert_{L^p(\mathbb{R}^n,\gamma(\mathbb{R}^n,\mathbb{B}))} &\le C \Vert G_{P,\alpha+1;\mathbb{B}}^H (T_\textrm{m} f) \Vert_{L^p(\mathbb{R}^n,\mathbb{B})}\\
    &\le C \int_{\mathbb{R}^n} \Vert \mathcal{M}_\alpha (\cdot,u) G_{P,1;\mathbb{B}}^H (L^{i \frac{u}{2}}f) \Vert_{L^p(\mathbb{R}^n,\gamma(\mathcal{H}^n,\mathbb{B}))} du\\
    &\le C \int_{\mathbb{R}^n} \sup_{t\in (0,\infty)^n} |\mathcal{M}_\alpha (\cdot,u)| \Vert  L^{i \frac{u}{2}} \Vert_{L^p(\mathbb{R}^n,\mathbb{B}) \to L^p(\mathbb{R}^n,\mathbb{B})} du \Vert f \Vert_{L^p(\mathbb{R}^n,\mathbb{B})}.
\end{align*}
We conclude that $T_\textrm{m}$ can be extended from $S(\mathbb{R}^n)\otimes \mathbb{B}$ to $L^p(\mathbb{R}^n,\mathbb{B})$ as a bounded operator from $L^p(\mathbb{R}^n,\mathbb{B})$ into itself.

\section{Proof of Theorem \ref{theorem1.6}}

According to \cite[Theorem 2.5.1]{TA} if $X$ is a UMD Banach space, $\alpha\in\mathbb{R}$ and $1<q<\infty$, the imaginary power $H^{i\alpha}$ of the Hermite operator is bounded from $L^q(\mathbb{R},X)$ into itself and, for every $\sigma \in\left(\frac{\pi}{2},\pi\right]$, there exists $C=C_{q,\sigma}>0$ such that
$$
\Vert H^{i\alpha} \Vert_{L^q(\mathbb{R},X) \to L^q(\mathbb{R},X)} \le C_{q,\sigma} e^{\sigma |\alpha|}.
$$
Moreover, if $\mathcal{H}$ is a Hilbert space then, for every $\alpha\in\mathbb{R}$, $H^{i\alpha}$ is bounded from $L^2(\mathbb{R},\mathcal{H})$ into itself and $\Vert H^{i\alpha} \Vert_{L^2(\mathbb{R},\mathcal{H}) \to L^2(\mathbb{R},\mathcal{H})} \le 1$.

Assume now that $\mathbb{B}$ is isomorphic to a closed subquotient of a complex interpolation space $[\mathcal{H},X]_\theta$, where $0<\theta<1$, $\mathcal{H}$ is a Hilbert space and $X$ is a UMD Banach space. Although the Hermite semigroup is not a diffusion semigroup we can proceed as in the proof of \cite[Corollary 2.5.3]{TA} to obtain that if $\left|\frac{2}{p}-1\right| < \theta$ and $0\le \Omega< \frac{\pi}{2}(1-\theta)$, there exists $\omega < \frac{\pi}{2} - \Omega$ such that
$$
\Vert H^{i\alpha} \Vert_{L^q(\mathbb{R},\mathbb{B}) \to L^q(\mathbb{R},\mathbb{B})} \le C e^{\omega |\alpha|},\;\; \alpha\in\mathbb{R},
$$
where $C$ does not depend on $\alpha$.

Hence, if $u=(u_1,\dots,u_n)\in\mathbb{R}^n$ we conclude that
\begin{equation}\label{n25}
\Vert L^{iu} \Vert_{L^q(\mathbb{R},\mathbb{B})} \le C e^{\omega \sum_{i=1}^n|u_i|},
\end{equation}
where $C$ does not depend on $u$.

By proceeding as in the proof of \cite[Theorem 3]{Me} and working coordinate to coordinate we can obtain that
$$
\sup_{t\in (0,\infty)^n} |\mathcal{M}_\gamma (t,u)| \le C \prod_{j=1}^n (1+|u_j|) |\Gamma (\gamma_j-iu_j)|,\;\;u=(u_1,\dots,u_n)\in\mathbb{R}^n.
$$
Then, if $\left|\frac{2}{p}-1\right|<\theta$ from (\ref{n25}) it follows that
$$
\int_{\mathbb{R}^n} \sup_{t\in (0,\infty)^n} |\mathcal{M}_\gamma (t,u)| \Vert L^{i \frac{u}{2}} \Vert_{L^p(\mathbb{R}^n,\mathbb{B}) \to L^p(\mathbb{R}^n,\mathbb{B})} du \le C \int_{\mathbb{R}^n} \prod_{j=1}^n (1+|u_j|) e^{-\left(\frac{\pi}{2}-\omega\right) \sum_{i=1}^n |u_i|} du <\infty.
$$
By using Theorem \ref{theorem1.5} we deduce that the multiplier operator $T_\textrm{m}$ is bounded from $L^p(\mathbb{R}^n,\mathbb{B})$ into itself provided that $\left|\frac{2}{p}-1\right|<\theta$.

\section{Proof of Theorem \ref{theorem1.7}}

Suppose that $\mathbb{B}$ is isomorphic to a closed subquotient of $[X,\mathcal{Q}]_\theta$ where $\theta \in (0,1)$, $X$ is a UMD Banach space, and $\mathcal{Q}$ is a Hilbert space, and that $\mathbb{B}$ has the property ($\alpha$). We are going to show that $\tilde W_{H,\ell}^p (\mathbb{R}^n,\mathbb{B}) = W_{H,\ell}^p (\mathbb{R}^n,\mathbb{B}) = L_{H,\ell}^p (\mathbb{R}^n,\mathbb{B})$, $\ell\in\mathbb{N}\setminus\{0\}$ and $1<p<\infty$, being $\left\vert \frac{2}{p}-1\right\vert <\theta$.

Let $1<p<\infty$ and $\ell\in\mathbb{N}$. By proceeding as in the proof of \cite[Proposition 1]{BT1} we can see that the linear space $\mathfrak{F}_H\otimes \mathbb{B}$ is dense in $\tilde W_{H,\ell}^p (\mathbb{R}^n,\mathbb{B})$, $W_{H,\ell}^p (\mathbb{R}^n,\mathbb{B})$ and $L_{H,\ell}^p (\mathbb{R}^n,\mathbb{B})$. Here by $\mathfrak{F}_H$ we denote the linear space $\textup{span}\{h_k\}_{k\in\mathbb{N}^n}$ generated by $\{h_k\}_{k\in\mathbb{N}^n}$.

Firstly we prove that $\tilde W_{H,\ell}^p (\mathbb{R}^n,\mathbb{B}) = L_{H,\ell}^p (\mathbb{R}^n,\mathbb{B})$. We consider the one dimensional situation and the first order Hermite-Riesz transform $\tilde R$ defined by
$$
\tilde R f= A_{-1} H^{-\frac12} f,\;\; f \in \mathfrak{F}_H.
$$
This operator can be extended to $L^2(\mathbb{R})$ as a bounded operator from $L^2(\mathbb{R})$ into itself. By denoting this extension again $\tilde R$ we have that
$$
\tilde R f = \sum_{k=0}^\infty \frac{\sqrt{2(k+1)}}{\sqrt{2k+1}} c_k(f) h_{k+1},\;\; f \in L^2(\mathbb{R}).
$$
By defining $\tilde R$ on $L^2(\mathbb{R})\otimes \mathbb{B}$ in a usual way, according to \cite[Theorem 2.3]{AT}, since $\mathbb{B}$ is a UMD Banach space, $\tilde R$ can be extended from $(L^2(\mathbb{R})\cap L^p(\mathbb{R}))  \otimes \mathbb{B}$ to $L^p(\mathbb{R},\mathbb{B})$ as a bounded operator from $L^p(\mathbb{R,\mathbb{B}})$ into itself.

We define the multiplier operator $T_\textrm{m}$ by
$$
T_\textrm{m} f = \sum_{k=0}^\infty \textrm{m}(\lambda_k) c_k(f) h_k,\;\; f \in L^2(\mathbb{R}),
$$
where $\textrm{m}(z) = \sqrt{\frac{z}{z+1}}$, $z\in \Gamma_{\frac{\pi}{2}}$. It is clear that $\textrm{m}$ is a bounded holomorphic function in $\Gamma_{\frac{\pi}{2}}$. Then, according to Theorem \ref{theorem1.6}, $T_\textrm{m}$ can be extended to $L^p(\mathbb{R,\mathbb{B}})$ as a bounded operator from $L^p(\mathbb{R,\mathbb{B}})$ into itself.

We denote by $S_+$ the forward shift operator defined by
$$
S_+ (f) = \sum_{k=0}^\infty c_k(f) h_{k+1},\;\; f \in L^2(\mathbb{R}).
$$
Since $S_+ = \tilde R T_\textrm{m}$, $S_+$ can be extended to $L^p(\mathbb{R,\mathbb{B}})$ as a bounded operator from $L^p(\mathbb{R,\mathbb{B}})$ into itself.

By $S_-$ we denote the backward shift operator defined by
$$
S_- (f) = \sum_{k=1}^\infty c_k(f) h_{k-1},\;\; f \in L^2(\mathbb{R}).
$$
According to \cite[Lemma 2.3]{GrafYang}, since $\mathbb{B}^*$ is a UMD Banach space with the property ($\alpha$), we get
\begin{align*}
    ||S_- (f) ||_{L^p(\mathbb{R}^n,\mathbb{B})} & = \sup_{\substack{g\in \mathfrak{F}_H \otimes \mathbb{B}^* \\ ||g||_{L^{p'}(\mathbb{R}^n,\mathbb{B}^*)} \le 1}} \left| \int_{\mathbb{R}} \langle g(y), S_- (f)(x) \rangle_{\mathbb{B}^*,\mathbb{B}} dx \right|\\
    &= \sup_{\substack{g\in \mathfrak{F}_H \otimes \mathbb{B}^* \\ ||g||_{L^{p'}(\mathbb{R}^n,\mathbb{B}^*)} \le 1}} \left| \int_{\mathbb{R}} \langle S_+ (g)(x) , f(x) \rangle_{\mathbb{B}^*,\mathbb{B}} dx \right|\\
    &\le \sup_{\substack{g\in \mathfrak{F}_H \otimes \mathbb{B}^* \\ ||g||_{L^{p'}(\mathbb{R}^n,\mathbb{B}^*)} \le 1}} \Vert S_+ (g) \Vert_{L^{p'}(\mathbb{R},\mathbb{B}^*)} \Vert f \Vert_{L^{p}(\mathbb{R},\mathbb{B})}\\
    &\le C \Vert f \Vert_{L^{p}(\mathbb{R},\mathbb{B})},\;\; f \in \mathfrak{F}_H\otimes \mathbb{B}.
\end{align*}
Hence, the operator $S_-$ can be extended to $L^{p}(\mathbb{R},\mathbb{B})$ as a bounded operator from $L^{p}(\mathbb{R},\mathbb{B})$ into itself.

We now come back to the $n$-dimensional situation. Let $m=(m_1,\dots,m_n)\in\mathbb{N}^n$ and $j\in \mathbb{N}$, $1\le j \le n$. We consider the shift operator $S_{m,j}$ defined by
$$
S_{m,j} f = \sum_{\substack{k=(k_1,\dots,k_n)\in \mathbb{N}^n\\ k_l\ge m_l,\; l=1,\dots, j}} c_k(f) h_{k_1-m_1}(x_1)\cdots h_{k_j-m_j}(x_j) h_{k_{j+1}+m_{j+1}}(x_{j+1})\cdots h_{k_n+m_n}(x_n),
$$
where $f\in L^2(\mathbb{R}^n)$. By taking into account the $L^p$-boundedness properties for the one dimensional shifts we can deduce that the operator $S_{m,j}$ can be extended to $L^p(\mathbb{R}^n,\mathbb{B})$ as a bounded operator from $L^p(\mathbb{R}^n,\mathbb{B})$ into itself.

We define the Hermite Riesz transform $R_{m,j}$ as follows
$$
R_{m,j} f = A_1^{m_1}\cdots A_j^{m_j}A_{-(j+1)}^{m_{j+1}} \cdots A_{-n}^{m_n} H^{-\frac{|m|}{2}} f,\;\; f\in \mathfrak{F}_H.
$$

According to \cite[(8)]{BT2} we have that, for every $k=(k_1,\dots,k_n)\in \mathbb{N}^n$, $k_\ell\ge m_\ell$, $\ell=1,\dots,j$,
\begin{eqnarray*}
\lefteqn{A_1^{m_1}\cdots A_j^{m_j}A_{-(j+1)}^{m_{j+1}} \cdots A_{-n}^{m_n} H^{-\frac{|m|}{2}} h_k(x)} \\
&=& \frac{\prod_{\ell=1}^j(\prod_{s=0}^{m_\ell-1} \sqrt{2(k_\ell-s)}) \prod_{\ell=j+1}^n(\prod_{s=1}^{m_\ell} \sqrt{2(k_\ell+s)})}{(2|k|+n)^{\frac{|m|}{2}}}\\
&&\times h_{k_1-m_1,\dots,k_j-m_j,k_{j+1}+m_{j+1},\cdots, k_n+m_n}(x),\;\; x \in \mathbb{R}^n.
\end{eqnarray*}
Then, for every $f\in \mathfrak{F}_H$, by taking into account that the sum is finite, we can write
\begin{eqnarray*}
R_{m,j} f &=& \sum_{\substack{k=(k_1,\dots,k_n)\in \mathbb{N}^n\\ k_\ell\ge m_\ell,\; \ell=1,\dots, j}} c_k (f) \frac{\prod_{\ell=1}^j(\prod_{s=0}^{m_\ell-1} \sqrt{2(k_\ell-s)}) \prod_{\ell=j+1}^n(\prod_{s=1}^{m_\ell} \sqrt{2(k_\ell+s)})}{(2|k|+n)^{\frac{|m|}{2}}}\\
&&\hspace{3cm}\times h_{k_1-m_1,\dots,k_j-m_j,k_{j+1}+m_{j+1},\cdots, k_n+m_n}.
\end{eqnarray*}
We consider the multiplier operator $T_\textrm{m}$ defined by
\begin{eqnarray*}
T_\textrm{m} f &=& \sum_{k=(k_1,\dots,k_n)\in \mathbb{N}^n} c_k (f) \frac{\prod_{\ell=1}^j(\prod_{s=0}^{m_\ell-1} \sqrt{2(k_\ell+m_\ell-s)}) \prod_{\ell=j+1}^n(\prod_{s=1}^{m_\ell} \sqrt{2(k_\ell+s)})}{(2|k|+ 2\sum_{\ell=1}^j m_\ell +n)^{\frac{|m|}{2}}}\\
&&\hspace{3cm}\times h_k,\;\; f\in \mathfrak{F}_H.
\end{eqnarray*}
Here
$$
\textrm{m}(z_1,\dots,z_n) = \frac{\prod_{\ell=1}^j(\prod_{s=0}^{m_\ell-1} \sqrt{z_\ell+2(m_\ell-s)-1}) \prod_{\ell=j+1}^n(\prod_{s=1}^{m_\ell} \sqrt{z_\ell+2s-1})}{(z_1+\cdots+z_n+ 2\sum_{\ell=1}^j m_\ell)^{\frac{|m|}{2}}},
$$
where $(z_1,\dots,z_n) \in \Gamma_{\frac{\pi}{2}}$. Suppose that $j\ne 0$. Then, the function $\textrm{m}$ satisfies the hypothesis in Theorem \ref{theorem1.6} that assures that the operator $T_\textrm{m}$ can be extended to $L^p(\mathbb{R}^n,\mathbb{B})$ as a bounded operator from $L^p(\mathbb{R}^n,\mathbb{B})$ into itself. We can write
$$
R_{m,j} f = S_{\bar \beta^{n-j},j}\circ T_\textrm{m}\circ S_{\beta^j,j} f,\;\; f\in \mathfrak{F}_H,
$$
where $\beta^j=(m_1,\dots,m_j,0,\dots,0)$ and $\bar\beta^{n-j}=(0,\dots,0,m_{j+1},\dots,m_n)$.. Hence, $R_{m,j}$ can be extended to $L^p(\mathbb{R}^n,\mathbb{B})$ as a bounded operator from $L^p(\mathbb{R}^n,\mathbb{B})$ into itself.

Assume now that $j=0$. We consider the function
$$
\textrm{m}(z_1,\dots,z_n) = \frac{\prod_{\ell=1}^n \prod_{s=1}^{m_\ell} \sqrt{z_\ell+2s-\frac12}}{(z_1+\cdots+z_n+\frac{n}{2})^{\frac{|m|}{2}}},\;\;(z_1,\dots,z_n) \in \Gamma_{\frac{\pi}{2}},
$$
and the multivariate $(H-\frac12)$-multiplier operator
$$
T_\textrm{m}^{\frac12} (f) = \sum_{k=(k_1,\dots,k_n)\in \mathbb{N}^n} c_k(f) \frac{\prod_{\ell=1}^n \prod_{s=1}^{m_\ell} \sqrt{2(k_j+s)}}{(2|k|+n)^{\frac{|m|}{2}}} h_k,\;\; f\in \mathfrak{F}_H.
$$
According to Theorem \ref{theorem1.6} (see Remark \ref{rem}) the operator $T_\textrm{m}^{\frac12}$ can be extended to $L^p(\mathbb{R}^n,\mathbb{B})$ as a bounded operator from $L^p(\mathbb{R}^n,\mathbb{B})$ into itself. Since
$$
R_{m,0} f = S_{m,0} \circ T_\textrm{m}^{\frac12} (f),\;\; f\in \mathfrak{F}_H,
$$
$R_{m,0}$ can be extended to $L^p(\mathbb{R}^n,\mathbb{B})$ as a bounded operator from $L^p(\mathbb{R}^n,\mathbb{B})$ into itself.

The above arguments allow us to establish the following property.
\begin{propo}\label{proposition6.1}
Suppose that $\mathbb{B}$ is isomorphic to a closed subquotient of $[X,\mathcal{Q}]_\theta$, where $\theta\in(0,1)$, $X$ is a UMD Banach space and $\mathcal{Q}$ is a Hilbert space and that $\mathbb{B}$ has the property ($\alpha$) and $1<p<\infty$ being $\left| \frac{2}{p}-1 \right|<\theta$. If $m=(m_1,\dots,m_n)\in\mathbb{N}^n$ and $j=(j_1,\dots,j_n)\in\mathbb{Z}^n$, being $|j_i|=i$, $i=1,\dots,n$, the Hermite Riesz transform $R_m^j$ defined by
$$
R_m^j f = \prod_{i=1}^n A_{j_i}^{m_i} H^{-\frac{|m|}{2}} f,\;\; f\in \mathfrak{F}_H,
$$
where $|m|=\sum_{i=1}^n m_i$, can be extended to $L^p(\mathbb{R}^n,\mathbb{B})$ as a bounded operator from $L^p(\mathbb{R}^n,\mathbb{B})$ into itself.
\end{propo}

By taking into account \cite[Lemma 2.3]{GrafYang} we now can prove, by proceeding as in the proof of \cite[Theorem 4]{BT1}, that $\tilde W_{H,\ell}^p (\mathbb{R}^n,\mathbb{B}) = L_{H,\ell}^p (\mathbb{R}^n,\mathbb{B})$.

In order to see that $W_{H,\ell}^p (\mathbb{R}^n,\mathbb{B}) = L_{H,\ell}^p (\mathbb{R}^n,\mathbb{B})$ we can use the ideas presented in the proof of \cite[Theorem 1]{BT2}. To show that $L_{H,\ell}^p (\mathbb{R}^n,\mathbb{B}) \hookrightarrow W_{H,\ell}^p (\mathbb{R}^n,\mathbb{B})$ it is sufficient to replace the operator $A_j$ by $A_{-j}$, $1\le j \le n$, in \cite[p. 157]{BT2}.

If $j\in \mathbb{Z}$, $1\le |j| \le n$, we consider the Hermite Riesz transform $R_j^\ell$ by
$$
R_j^\ell f = A_j^\ell H^{-\frac{\ell}{2}} f,\;\; f\in \mathfrak{F}_H,
$$
and we define the operator $\tau_\ell= \sum_{j=1}^n R_j^\ell R_{-j}^\ell$. According to Proposition \ref{proposition6.1} the operator $\tau_\ell$ is bounded from $L^p(\mathbb{R}^n,\mathbb{B})$ into itself. By using \cite[Lemma 4 and (8)]{BT2} we obtain
\begin{align*}
    \tau_\ell f &= \sum_{j=1}^n A_j^\ell H^{-\frac{\ell}{2}} A_{-j}^\ell H^{-\frac{\ell}{2}} f \\
    &= (H+2)^{-\frac{\ell}{2}} \sum_{j=1}^n A_j^\ell A_{-j}^\ell H^{-\frac{\ell}{2}} f \\
    &= \sum_{j=1}^n \sum_{l=(l_1,\dots,l_n)\in \mathbb{N}^n} 2^\ell \frac{\prod_{m=1}^\ell (l_j+m)}{(2|l|+n+2)^{\frac{\ell}{2}} (2|l|+n)^{\frac{\ell}{2}}} c_l(f) h_l\\
    &= \sum_{l=(l_1,\dots,l_n)\in \mathbb{N}^n} \left( 2^\ell \frac{\sum_{j=1}^n \prod_{m=1}^\ell (l_j+m)}{(2|l|+n+2)^{\frac{\ell}{2}} (2|l|+n)^{\frac{\ell}{2}}} \right) c_l(f) h_l,\;\; f \in \mathfrak{F}_H
\end{align*}
We consider the function
$$
\textrm{m}(z_1,\dots,z_n) = \frac{(z_1+\cdots+z_n+2)^{\frac{\ell}{2}}{(z_1+\cdots+z_n)^{\frac{\ell}{2}}}}{\sum_{j=1}^n \prod_{m=1}^\ell (\frac{z_j}{2}+m-\frac12)},\;\; z=(z_1,\dots,z_n) \in \Gamma_{\frac{\pi}{2}}.
$$
According to Theorem \ref{theorem1.6}, the multiplier $T_\textrm{m}$ associated with $\textrm{m}$ is bounded from $L^p(\mathbb{R}^n,\mathbb{B})$ into itself.

Hence, we obtain
$$
\Vert g \Vert_{L^p(\mathbb{R}^n,\mathbb{B})} \le \Vert T_\textrm{m} \tau_\ell g \Vert_{L^p(\mathbb{R}^n,\mathbb{B})} \le C \sum_{j=1}^n \Vert A_{-j}^\ell H^{-\frac{\ell}{2}}g \Vert_{L^p(\mathbb{R}^n,\mathbb{B})},\;\; g \in \mathfrak{F}_H.
$$
Then,
$$
\Vert f \Vert_{L_{H,\ell}^p (\mathbb{R}^n,\mathbb{B})} \le C \Vert f \Vert_{W_{H,\ell}^p (\mathbb{R}^n,\mathbb{B})},\;\; f \in \mathfrak{F}_H.
$$
We conclude that $W_{H,\ell}^p (\mathbb{R}^n,\mathbb{B})$ is continuously contained in $L_{H,\ell}^p (\mathbb{R}^n,\mathbb{B})$.

Thus, the proof of this theorem is finished.

\section{Proof of Theorem \ref{nteo8}}

In order to show this result we can use the ideas developed in \cite{BFRTT}. We consider the space $\mathfrak{F}_H = \textup{span}\;\{h_k\}_{k\in \mathbb{N}^n}$ generated by $\{h_k\}_{k\in \mathbb{N}^n}$. Also we define, for every $\beta>0$ and $k\in\mathbb{N}$ such that $k>\beta$, the operators
\begin{align*}
G_{P,\beta;\mathbb{B}}^H (f)(t,x) & =  t^\beta \partial_t^\beta P_t^H(f)(x),\;\; t\in (0,\infty)~\textup{and}~x\in\mathbb{R}^n,\\
\intertext{and}
G_{P,\beta,k;\mathbb{B}}^H (f)(t,x) & =  t^{k-\beta} \partial_t^k P_t^H(f)(x),\;\; t\in (0,\infty)~\textup{and}~x\in\mathbb{R}^n,
\end{align*}
for every $f\in \mathfrak{F}_H\otimes \mathbb{B}$.

Let $1<p<\infty$, $\beta>0$ and $k\in\mathbb{N}$ such that $k>\beta$. According to \cite[Lemma 2.2]{BFRTT} we have that
$$
G_{P,k-\beta;\mathbb{B}}^H (f) = G_{P,\beta,k;\mathbb{B}}^H (H^{-\frac{\beta}{2}} f),\;\; f \in \mathfrak{F}_H\otimes \mathbb{B}.
$$
Then, from \cite[Theorem 1]{BCCFR1} it follows that
\begin{equation}\label{n7.1}
\frac{1}{C} \Vert f \Vert_{L_{H,\frac{\beta}{2}}^p(\mathbb{R}^n,\mathbb{B})} \le \Vert G_{P,\beta,k;\mathbb{B}}^H (f) \Vert_{L^p(\mathbb{R}^n,\gamma(\mathcal{H}^1,\mathbb{B}))} \le C \Vert f \Vert_{L_{H,\frac{\beta}{2}}^p(\mathbb{R}^n,\mathbb{B})},\;\; f\in \mathfrak{F}_H\otimes \mathbb{B}.
\end{equation}

The space $\mathfrak{F}_H\otimes \mathbb{B}$ is dense in $L_{H,\frac{\beta}{2}}^p(\mathbb{R}^n,\mathbb{B})$. Indeed, let $f\in L^p(\mathbb{R}^n,\mathbb{B})$. Since the operator $H^{-\frac{\beta}{2}}$ is bounded from $L^p(\mathbb{R}^n,\mathbb{B})$ into itself and $\mathfrak{F}_H\otimes \mathbb{B}$ is dense in $L^p(\mathbb{R}^n,\mathbb{B})$, there exists a sequence $\{g_k\}_{k=1}^\infty \subset \mathfrak{F}_H\otimes \mathbb{B}$ such that $g_k \to H^{-\frac{\beta}{2}} f$, as $k\to \infty$, in $L^p(\mathbb{R}^n,\mathbb{B})$. We define, for every $k\in\mathbb{N}\setminus\{0\}$,
$$
\Lambda_k = \sum_{l\in \mathbb{N}^n} c_l(g_k) h_l (2l+|n|)^{\frac{\beta}{2}}.
$$
Note that $\Lambda_k \in \mathfrak{F}_H$ and $H^{-\frac{\beta}{2}}(\Lambda_k) = g_k$, $k\in \mathbb{N}\setminus\{0\}$. We have that $\Lambda_k \to f$, as $k\to \infty$, in $L^p_{H,\frac{\beta}{2}}(\mathbb{R}^n,\mathbb{B})$.

Then, (\ref{n7.1}) implies that the operator $G_{P,\beta,k;\mathbb{B}}^H$ can be extended to $L^p_{H,\frac{\beta}{2}}(\mathbb{R}^n,\mathbb{B})$ as a bounded operator $\tilde{G}_{P,\beta,k;\mathbb{B}}^H$ from $L^p_{H,\frac{\beta}{2}}(\mathbb{R}^n,\mathbb{B})$ into $L^p(\mathbb{R}^n,\gamma(\mathcal{H}^1,\mathbb{B}))$.

By proceeding as in the proof of Theorem \ref{theorem1.4} we can see that, for every $f\in L^p_{H,\frac{\beta}{2}}(\mathbb{R}^n,\mathbb{B})$, $\tilde{G}_{P,\beta,k;\mathbb{B}}^H (f)(x) = G_{P,\beta,k;\mathbb{B}}^H (f)(\cdot,x)$, a.e. $x\in\mathbb{R}^n$. Then,
\begin{equation}\label{n14}
\frac{1}{C} \Vert f \Vert_{L^p_{H,\frac{\beta}{2}}(\mathbb{R}^n,\mathbb{B})} \le \Vert G_{P,\beta,k;\mathbb{B}}^H (f) \Vert_{L^p(\mathbb{R}^n,\gamma(\mathcal{H}^1,\mathbb{B}))} \le C \Vert f \Vert_{L^p_{H,\frac{\beta}{2}}(\mathbb{R}^n,\mathbb{B})},
\end{equation}
for every $f\in L^p_{H,\frac{\beta}{2}}(\mathbb{R}^n,\mathbb{B})$.

Suppose now that $f\in L^p(\mathbb{R}^n,\mathbb{B})$ and $G_{P,\beta,k;\mathbb{B}}^H (f)\in L^p(\mathbb{R}^n,\gamma(\mathcal{H}^1,\mathbb{B}))$. We are going to see that $f\in L^p_{H,\frac{\beta}{2}}(\mathbb{R}^n,\mathbb{B})$. Let $\delta >0$. We define
$$
F_\delta = \sum_{k \in \mathbb{N}^n} (2|k|+n)^{\frac{\beta}{2}} e^{-\delta(2|k|+n)^{\frac12}} c_k^H(f) h_k.
$$
By \cite[(2.2)]{StemTo1}, $F_\delta \in L^p(\mathbb{R}^n,\mathbb{B})$. Since $H^{-\frac{\beta}{2}} F_\delta=P_\delta^H (f) \in L^p(\mathbb{R}^n,\mathcal{H})$, $P_\delta^H (f) \in L^p_{H,\frac{\beta}{2}}(\mathbb{R}^n,\mathbb{B})$ and $\Vert P_\delta^H(f) \Vert_{L^p_{H,\frac{\beta}{2}}(\mathbb{R}^n,\mathbb{B})} = \Vert F_\delta \Vert_{L^p(\mathbb{R}^n,\mathbb{B})}$. Hence, (\ref{n14}) implies that
$$
\Vert F_\delta \Vert_{L^p(\mathbb{R}^n,\mathbb{B})} \le C \Vert G_{P,\beta,k;\mathbb{B}}^H (P_\delta^H(f)) \Vert_{L^p(\mathbb{R}^n,\gamma(\mathcal{H}^1,\mathbb{B}))}.
$$
Also, we can write
\begin{align*}
G_{P,\beta,k;\mathbb{B}}^H (P_\delta^H(f)) (x) &= t^{k-\beta} \partial_t^k P_t^H (P_\delta^H(f))(x) = t^{k-\beta} \partial_t^k P_\delta^H(P_t^H (f))(x)\\
    &= P_\delta^H (t^{k-\beta} \partial_t^k P_t^H (f))(x),\;\; x \in \mathbb{R}^n~\textup{and}~t\in (0,\infty),
\end{align*}
and, then
$$
\Vert G_{P,\beta,k;\mathbb{B}}^H (P_\delta^H(f)) (x) \Vert_{\gamma(\mathcal{H}^1,\mathbb{B})} \le P_\delta^H \left(\Vert t^{k-\beta} \partial_t^k P_t^H (f) \Vert_{\gamma(\mathcal{H}^1,\mathbb{B})}\right)(x),\;\; \textup{a.e.}~x\in\mathbb{R}^n.
$$
Since $P_\delta^H$ is contractive in $L^p(\mathbb{R}^n)$ we obtain
$$
\Vert G_{P,\beta,k;\mathbb{B}}^H (P_\delta^H(f)) \Vert_{L^p(\mathbb{R}^n,\gamma(\mathcal{H}^1,\mathbb{B}))} \le \Vert G_{P,\beta,k;\mathbb{B}}^H (f) \Vert_{L^p(\mathbb{R}^n,\gamma(\mathcal{H}^1,\mathbb{B}))}.
$$
Hence, the set $\{F_\delta\}_{\delta>0}$ is bounded in $L^p(\mathbb{R}^n,\mathbb{B})$.

Since $\mathbb{B}$ is a UMD Banach space, the dual of $L^p(\mathbb{R}^n,\mathbb{B})$ is $L^{p'}(\mathbb{R}^n,\mathbb{B}^*)$, where $p'=\frac{p}{p-1}$. By Banach-Alaoglu's Theorem there exists a decreasing sequence $\{\delta_k\}_{k\in \mathbb{N}\setminus\{0\}} \subset (0,\infty)$ and $F\in L^p(\mathbb{R}^n,\mathbb{B})$ such that $\delta_k \to 0$, as $k\to \infty$,
$$
\int_{\mathbb{R}^n} \langle g(x), F_{\delta_k}(x)\rangle_{\mathbb{B}^*,\mathbb{B}} dx \to \int_{\mathbb{R}^n} \langle g(x), F(x)\rangle_{\mathbb{B}^*,\mathbb{B}} dx,\;\;\textup{as}~k\to \infty,
$$
for every $g\in L^{p'}(\mathbb{R}^n,\mathbb{B}^*)$, and $\Vert F \Vert_{L^p(\mathbb{R}^n,\mathbb{B})} \le C \Vert G_{P,\beta,k;\mathbb{B}}^H (f) \Vert_{L^p(\mathbb{R}^n,\gamma(\mathcal{H}^1,\mathbb{B}))}$. By taking into account that $H^{-\frac{\beta}{2}}$ is a bounded operator from $L^{p}(\mathbb{R}^n,\mathbb{B})$ into itself, we get, for every $g\in L^{p'}(\mathbb{R}^n,\mathbb{B}^*)$,
$$
\int_{\mathbb{R}^n} \langle g(x), P_{\delta_k}^H(f)(x)\rangle_{\mathbb{B}^*,\mathbb{B}} dx \to \int_{\mathbb{R}^n} \langle g(x), H^{-\frac{\beta}{2}}(F)(x)\rangle_{\mathbb{B}^*,\mathbb{B}} dx,\;\;\textup{as}~k\to \infty.
$$
Since $P_{\delta_k}^H(f) \to f$, as $k\to \infty$, in $L^p(\mathbb{R}^n,\mathbb{B})$, we conclude that $f= H^{-\frac{\beta}{2}}(F)$. Hence, $f\in L^p_{H,\frac{\beta}{2}}(\mathbb{R}^n,\mathbb{B})$ and
$$
\Vert f \Vert_{L^p_{H,\frac{\beta}{2}}(\mathbb{R}^n,\mathbb{B})} \le C \Vert G_{P,\beta,k;\mathbb{B}}^H (f) \Vert_{L^p(\mathbb{R}^n,\gamma(\mathcal{H}^1,\mathbb{B}))}.
$$
Thus, the proof of this theorem is finished.

\begin{remark}
Classical Sobolev and potential spaces in a Banach valued setting have been studied, for instance, in \cite{Buk}, \cite{PW} and \cite{SS}. If $\mathbb{B}$ is a UMD Banach space the classical Sobolev space $W_\ell^p (\mathbb{R}^n,\mathbb{B})$ coincides with the classical potential space $L^p_\ell (\mathbb{R}^n,\mathbb{B})$, for every $1<p<\infty$ and $\ell \in \mathbb{N}$. In this case it is not necessary to consider UMD spaces with the property ($\alpha$) as in our Hermite setting. As it is written in \cite[Theorem 3, p. 135]{Stein} the proof of the equality $W^p_\ell (\mathbb{R}^n,\mathbb{B}) = L^p_\ell (\mathbb{R}^n,\mathbb{B})$, $1<p<\infty$ and $\ell \in \mathbb{N}$, when $\mathbb{B}=\mathbb{C}$ (in the scalar case) also works when $\mathbb{B}$ is a UMD Banach space. Indeed, it is sufficient to use Fourier multipliers in UMD Banach spaces (\cite{McConnell} and \cite{We}) and to take into account that some properties for the classical Riesz transforms hold true in UMD Banach valued settings. As Bongioanni and Torrea \cite{BT1} show in the Hermite case in order to show the equality $W^p_\ell (\mathbb{R}^n,\mathbb{C}) = L^p_\ell (\mathbb{R}^n,\mathbb{C})$ (even in the scalar case), multivariate Hermite multipliers are required. This fact leads us to consider UMD Banach spaces with the property ($\alpha$) in our Hermite case. It is an interesting question if the $(\alpha)$-property for $\mathbb B$ is necessary in order that the equality $W_{\ell}^p(\mathbb R^n,\mathbb C)=L^p_\ell(\mathbb R^n,\mathbb B)$, $1<p<\infty$ and $\ell \in \mathbb N$, holds.
\end{remark}


\end{document}